\input amstex
\documentstyle{amsppt}
\NoRunningHeads
\pagewidth{14cm}
\pageheight{19.1cm}
\magnification=1200
\loadmsbm
\loadbold

\define\End{\operatorname {End}}
\define\red{{\operatorname {red}}}
\define\rad{\operatorname {rad}}

\define\id{\operatorname {id}}
\define\ad{\operatorname {ad}}
\define\tr{\operatorname {tr}}

\define\spa{\operatorname {span}}
\define\bs{\boldsymbol}

\define\ol{\overline}
\define\vvp{\varphi}

\define\bk{\boldkey}

\define\der{\operatorname {Der}}

\NoBlackBoxes

\dedicatory Dedicated to Stephen Berman on the occasion of his 60th birthday
\enddedicatory

\author Jun Morita
\footnote{\leftline{Partially supported by a MONKASHO KAKENHI (2003--2004).}}
\endauthor

\affil Institute of Mathematics \\ University of Tsukuba \\ Tsukuba, Ibaraki, 305-8571 Japan \\
morita\@math.tsukuba.ac.jp\\
\\ Yoji Yoshii
\\ Department of Mathematics\\ North Dakota State University \\  Fargo, ND, 58105-5075  USA \\
yoji.yoshii\@ndsu.nodak.edu
\endaffil

\topmatter
\title Locally Extended Affine Lie Algebras
\endtitle

\abstract We propose a new simplified definition of extended affine Lie algebras (EALAs for short), and
also discuss a general version of extended affine Lie algebras, called locally extended affine Lie
algebras (LEALAs for short). We prove a conjecture by V. Kac for LEALAs. It turns out that the root
system of a LEALA becomes a locally finite version of extended affine root systems. Finally, several
examples of new EALAs and LEALAs are introduced, and we classify LEALAs of nullity 0 with connection to
locally finite split simple Lie algebras.
\endabstract

\endtopmatter

\document

\head Introduction
\endhead

Extended affine Lie algebras, or EALAs for short, were first introduced by H\o egh-Krohn and Torresani
in 1990 [HT] (under the name of irreducible quasi-simple Lie algebras), as a generalization of 
finite-dimensional simple Lie algebras and affine Kac-Moody Lie algebras over the complex numbers
$\Bbb C$. EALAs were systematically studied by Allison, Azam, Berman, Gao and Pianzola in [A-P]. They
proved the so-called Kac conjecture, 
which implies that the root systems of EALAs 
are examples of extended
affine root systems which were previously introduced by Saito [S] in 1985.

A natural question about the definition of original EALAs 
is the necessity of working over $\Bbb C$.   A
recent announcement by Neher [N2] has fixed this problem, taking as our base field an arbitrary field
$F$ of characteristic
$0$, and he has reported broad results about EALAs over $F$.

Our initial purpose was to give a more general (and simpler) definition of EALAs over $F$ and prove the
Kac conjecture. Roughly speaking, the conjecture says that a naturally induced  symmetric bilinear form
on the
$\Bbb Q$-span of the roots is positive semidefinite. It turns out that  the conjecture is true, and the
proof can be somewhat simplified with an argument in [AKY] (and of course with the original methods
in [A-P]). Thus one can define the {\it nullity} of an EALA as the $\Bbb Q$-dimension of the radical of
the form, as in the original theory. However,  in this note we do not need any assumptions 
for the set of roots $R$,
except for the irreducibility of the set of anisotropic roots. The original theory assumed
$R$ to be discrete, but there is no such concept in our setting because our base field is more general.
Also, an isotropic root can be isolated (e.g. all the roots can be isotropic in our setting). On the
other hand, Neher assumed that the additive group generated by isotropic roots has finite rank [N2]. As
a result, his EALAs are the tame EALAs of finite null rank in our sense.

We will give examples of an EALA of finite nullity but not having finite rank, and an EALA of infinite
nullity. The latter algebra is constructed in a way analogous to the construction of affine Kac-Moody
Lie algebras from an infinite-loop algebra
$$\frak g\otimes_F F [t_i^{\pm 1}]_{i\in\Bbb N},$$ where $\frak g$ is a finite-dimensional split simple
Lie algebra over $F$ and $F[t_i^{\pm 1}]_{i\in\Bbb N}$ is the ring of Laurent polynomials in 
infinitely many variables.  (One needs $[F:\Bbb Q]=\infty$ to make this loop algebra an EALA.) 
So it seems
reasonable to call the algebra an {\it EALA of nullity $\infty$} (see \S 5).

We have also noticed that the finite dimensionality of  our Cartan subalgebra 
is not used much for the
theory. So we exclude this assumption, and instead an axiom, (A3) in \S 1, is added. This roughly
says that $R$ is small enough to be captured by our nondegenerate form. We call this new algebra a {\it
locally extended affine Lie algebra} or a {\it LEALA} for short. We will prove the Kac conjecture for
LEALAs. An interesting phenomenon is that the roots of an EALA consist of a finite irreducible root
system and isotropic roots, while the roots of a LEALA consist of a locally finite irreducible root
system and isotropic roots. Thus a generalization of 
Saito's extended affine root systems
naturally comes up.

The so-called affine Lie algebras of infinite rank (see [K, \S 7.11]) or locally finite split simple
Lie algebras ([Stu] or [NS]) are our LEALAs of nullity 0.  (We expect that LEALAs will be an interesting
topic in the context of locally finite algebras and the recent work 
on locally finite root systems by Loos and Neher [LN].)

Also, in our axioms of a LEALA (or an EALA), even if there is no anisotropic root, 
the Kac conjecture still
holds. Namely the form becomes identically zero. 
We call such an algebra a {\it null system}. Heisenberg Lie
algebras with derivations are such examples (Example 7.1). Also, we construct an interesting null
system of nullity $\infty$ from a generalized Witt algebra and its dual module (Example 7.3). This null
system of finite null rank coincides with  a subalgebra (null part) of an EALA of maximal type
constructed in [BGK].

We classify LEALAs of nullity 0. The EALAs of nullity 0 are exactly the split central extensions of
finite-dimensional split simple Lie algebras. However, a new phenomenon comes into play for LEALAs.
We first show that the core of a LEALA $L$ of nullity 0 is a locally finite split simple Lie algebra,
and the centre of $L$ is always split. So a problem is to classify the centreless LEALAs of nullity 0
(while the centreless EALAs are exactly the finite-dimensional split simple Lie algebras). 
Locally finite split
simple Lie algebras do not exhaust the class. One of the structural difference comes from the fact
that  infinite-dimensional locally finite split simple Lie algebras have outer derivations.
Consequently, a centreless LEALA of nullity 0 is isomorphic to the semidirect product of a  locally
finite split simple Lie algebra with a certain family of outer derivations.

In the final section,  we discuss the relation between indecomposability and tameness of LEALAs.

Finally, we would like to thank Saied Azam, Karl-Hermann Neeb, and  Erhard Neher for several profound
insights and suggestions.

\head
\S 1 Definition of a LEALA
\endhead

Throughout the paper, let $F$ be a field of characteristic $0$. We say that a Lie algebra $L$ has a
{\it root decomposition} with respect to an abelian subalgebra $H$ if
$$L=\bigoplus_{\xi\in H^*}\ L_\xi,$$ where
$H^*$ is the dual space of $H$ and
$$L_\xi=\{x\in L\mid [h,x]=\xi(h)x\
\text{for all}\ h\in H\}.$$ An element of the set
$$R=R(H)=\{\xi\in H^*\mid L_\xi\neq 0\}$$ is called a {\it root}. We consider a category
$\Cal L$ of Lie algebras with such an $H$ and a symmetric invariant bilinear form ${\Cal B}$.  For
$(L,H,{\Cal B}), (L',H',{\Cal B}')\in\Cal L$, a morphism $\vvp:(L,H,{\Cal
B})\longrightarrow(L',H',{\Cal B}')$ is a Lie algebra homomorphism $\vvp:L\longrightarrow L'$ such that
$\vvp(H)\subset H'$ and ${\Cal B}(x,y)={\Cal B}'(\vvp(x),\vvp(y))$ for all $x,y\in L$.  An isomorphism
$\vvp:(L,H,{\Cal B})\longrightarrow(L',H',{\Cal B}')$ is a morphism such that $\vvp$ is a Lie algebra
isomorphism
$\vvp:L\longrightarrow L'$ with $\vvp(H)= H'$. If we take $(L_1,H_1,{\Cal B}_1), (L_2,H_2,{\Cal
B}_2)\in {\Cal L}$, then we have $$ (L_1,H_1,{\Cal B}_1) \oplus (L_2,H_2,{\Cal B}_2): =(L_1\oplus
L_2,H_1\oplus H_2,{\Cal B}_1\perp {\Cal B}_2)\in {\Cal L}.$$  If a triple $(L,H,{\Cal B})$ can never be
expressed as $(L_1,H_1,{\Cal B}_1) \oplus (L_2,H_2,{\Cal B}_2)$, a direct sum of nontrivial triples,
then $(L,H,{\Cal B})$ is called {\it indecomposable}.

We note that the Jacobi identity implies that
$[L_\xi,L_\eta]\subset L_{\xi+\eta}$, and the invariance of ${\Cal B}$ implies
${\Cal B}(L_\xi,L_\eta)=0$ unless $\xi+\eta=0$ for
$\xi,\eta\in H^*$.

A triple $(L,H,{\Cal B})\in\Cal L$ (or simply $L$) is called  {\it admissible} if it satisfies
\roster
\item"(A1)" $H$ is self-centralizing, i.e., $L_0=H$;

\item"(A2)" ${\Cal B}$ is nondegenerate;

Hence
$L_\xi$ and $L_{-\xi}$ form a nondegenerate pair relative to the form ${\Cal B}$, and ${\Cal
B}\mid_{H\times H}$ is nondegenerate. In particular, $-R=R$.

\item"(A3)" $R\subset H_{\Cal B}^*$, where
$H_{\Cal B}^*$ is the image of the canonical map
$H\longrightarrow H^*$ defined by $h\mapsto {\Cal B}(h,\cdot)$;

Hence there is the induced form on $H_{\Cal B}^*$ from ${\Cal B}$, simply denoted $(\cdot,\cdot)$. That
is, since
${\Cal B}\mid_{H\times H}$ is nondegenerate, we define
$(\xi,\eta):={\Cal B}(t_\xi,t_\eta)$ for $\xi,\eta\in H_{\Cal B}^*$, where $t_\xi$ (or $t_\eta$
similarly) is the unique element so that
$\xi(h)={\Cal B}(t_\xi,h)$ for all $h\in H$. Then, we have
$$[x,y]={\Cal B}(x,y)t_\xi
\tag$1.0$
$$ for $\xi\in R$, $x\in L_\xi$ and
$y\in L_{-\xi}$.

\endroster

\example{Remark 1.1}  (1) Under the assumptions on $H$, the normalizer and the centralizer of
$H$ coincides, and so  (A1) is equivalent to saying that $H$ is a Cartan subalgebra (for a more general
definition of Cartan subalgebras, see [NP]).  In fact, a Lie algebra $L=(L,H)$ satisfying (A1) is
called a {\it split Lie algebra} and
$H$ is called a {\it split Cartan subalgebra} of $L$ in [NS, Def.II.1].

(2) If $H$ is finite-dimensional, then (A3) is automatically satisfied, i.e., $H_{\Cal B}^*=H^*$.
\endexample

\enskip

We call an element of the set
$$ R^\times=R^\times(H,{\Cal B})=\{\xi\in R\mid (\xi,\xi)\neq 0\}
\enskip\text{(resp.}\enskip R^0=R^0(H,{\Cal B})=\{\xi\in R\mid (\xi,\xi)= 0\})
$$ an {\it anisotropic root} (resp.  an {\it isotropic root}). Note that $R$ or $R^\times$ can be
empty. For each nonzero $\xi\in R$, one chooses and fixes
$x_\xi\in L_\xi$ and
$x_{-\xi}\in L_{-\xi}$ such that
$[x_\xi,x_{-\xi}]=t_\xi$ by ($1.0$). If $\xi\in R^\times$, the Lie algebra generated by $x_\xi$ and
$x_{-\xi}$ is isomorphic to $sl_2(F)$. Letting $y_\xi:=\frac{2}{(\xi,\xi)}x_{-\xi}$ and
$h_\xi:=\frac{2}{(\xi,\xi)}t_\xi$, we call
$(x_\xi,h_\xi,y_\xi)$ an {\it $sl_2$-triplet for $\xi\in R^\times$}. Also, if $0\neq\xi\in R^0$, the
Lie algebra generated by $x_\xi$ and $x_{-\xi}$ is a $3$-dimensional Heisenberg Lie algebra. Letting
$y_\xi:=x_{-\xi}$, we call
$(x_\xi,t_\xi,y_{\xi})$ a  {\it Heisenberg-triplet for $0\neq\xi\in R^0$}. Thus, an admissible triple
is generated by a bunch of copies of $sl_2(F)$,
$3$-dimensional Heisenberg Lie algebras, and $H$.

Any symmetrizable Kac-Moody Lie algebra is an example of admissible triples with finite-dimensional
$H$.

\definition{Definition 1.2} An admissible triple $(L,H,{\Cal B})\in\Cal L$ is called  a {\it locally
extended affine Lie algebra} or a {\it LEALA} for short if it satisfies
\roster
\item"(A4)" $\ad x\in\End_F L$ is locally nilpotent for all $\xi\in R^\times$ and all $x\in L_\xi$, 

\item"(A5)" $R^\times$ is irreducible, i.e., $R^\times=R_1\cup R_2$ and $(R_1,R_2)=0$ imply
$R_1=\emptyset$ or
$R_2=\emptyset$.
\endroster 

Also:
\roster
\item"(i)" If $H$ is finite-dimensional, then $L$ is called an {\it extended affine Lie algebra}, an
{\it EALA} for short. 

\item"(ii)" If $R^\times=\emptyset$, then
$(L,H,{\Cal B})$ is called a {\it null LEALA} (or a {\it null EALA}) or simply a {\it null system}.

\item"(iii)" If the triple $(L,H,{\Cal B})$ is indecomposable, then it is called an {\it indecomposable
EALA}, an {\it indecomposable LEALA}   or an {\it indecomposable null system}.
 
\endroster
\enddefinition

\example{Remark 1.3}  Our EALAs generalize the original EALAs in [A-P].  For example, we only assume
the base field $F$ to be of characteristic 0, while  the base field in [A-P] is $\Bbb C$ (the complex
numbers). Also, even over
$\Bbb C$, we will show in \S 5 that our class of EALAs is wider than the original one, but we would
like to keep the same name EALA since the essential structure is same.  EALAs defined by Neher in [N2]
are tame EALAs of finite null rank in our sense (see Definition 3.11 and 9.1).
\endexample

\head
\S 2 Basic Properties
\endhead

Let $(L,H,{\Cal B})\in\Cal L$ be a LEALA. Using $sl_2$-theory, one can prove the following lemma in the
same way as in [A-P, Thm I.1.29]. (We do not need to assume (A5) in this section.)

\proclaim{Lemma 2.1} Let $\alpha\in R^\times$. Then:
\roster
\item For $\xi\in R$ we have
$\xi(h_\alpha) =
\frac{2(\xi,\alpha)}{(\alpha,\alpha)}\in\Bbb Z$.
\item The reflection
$\sigma_\alpha$ defined by
$\sigma_\alpha(\mu) = \mu - \mu(h_\alpha)\alpha$ for all $\mu \in H^\ast$ preserves $R$, that is,
$\sigma_\alpha(R) = R$. Also, $\sigma_\alpha$ is in the orthogonal group of the form on $H_{\Cal B}^*$.
\item If $k\in F$ and $k\alpha\in R$ then $k=0,\pm 1$.
\item
$\dim_FL_\alpha=1$.
\item For any $\xi\in R$ there exist two non-negative integers $d$, $u$ such that for any $n\in\Bbb Z$
we have
$\xi+n\alpha\in R$ if and only if
$-d\leq n\leq u$, where
$\xi - d\alpha,\ldots,\xi,\ldots,\xi + u\alpha$ is called the
$\alpha$-string through $\xi$. Moreover, $d-u=\xi(h_\alpha)=\frac{2(\xi,\alpha)}{(\alpha,\alpha)}$.
\endroster
\endproclaim

Also, we will use the following.

\proclaim{Lemma 2.2} Let $\alpha,\beta\in R^\times$. If $\alpha+\beta$ and $\alpha-\beta\in R$, then
$(\alpha,\alpha)=(\beta,\beta)$.
\endproclaim
\demo{Proof} Let $(x_{\alpha+\beta},t_{\alpha+\beta},x_{-\alpha-\beta})$ 
be a triplet satisfying $\Cal B(x_{\alpha+\beta},x_{-\alpha-\beta}) = 1$ (see (1.0)). 
Let $0\neq x\in L_{\alpha-\beta}$. Since $2\alpha,-2\beta\notin R$, we have
$[x_{\alpha+\beta},x]=0=[x_{-\alpha-\beta},x]$, and so
$((\alpha,\alpha)-(\beta,\beta))x = (\alpha-\beta,\alpha+\beta)x =(\alpha-\beta)(t_{\alpha+\beta})x =
[t_{\alpha+\beta},x] =[[x_{\alpha+\beta},x_{-\alpha-\beta}],x]
=[x_{\alpha+\beta},[x_{-\alpha-\beta},x]]- [x_{-\alpha-\beta},[x_{\alpha+\beta},x]] = 0$. Hence,
$(\alpha,\alpha)-(\beta,\beta)=0$.
\qed
\enddemo

\head
\S 3 Kac conjecture
\endhead

Let $(L,H,{\Cal B}) \in \Cal L$ be a LEALA. Here we scale the form ${\Cal B}$ on $L$, which induces a
scaling of the form
$(\cdot,\cdot)$. Namely, we put ${\Cal B}' = u {\Cal B}$ for some nonzero element $u
\in F$. This new form ${\Cal B}'$ induces another form on $H_{\Cal B}^*=H_{{\Cal B}'}^*$, called
$(\cdot,\cdot)'$. If $\alpha \in R$ with $\alpha \not= 0$, then
${\Cal B}(t_\alpha,h) = \alpha(h) = {\Cal B}'(t_\alpha',h) = u {\Cal B}(t_\alpha',h) = {\Cal B}(u
t_\alpha',h)$ and hence $t_\alpha = u t_\alpha'$, which implies
$(\alpha,\alpha)' = {\Cal B}'(t_\alpha',t_\alpha') = u {\Cal B}(t_\alpha',t_\alpha') = {\Cal B}(u
t_\alpha',t_\alpha') = {\Cal B}(t_\alpha,t_\alpha/u)  = {\Cal B}(t_\alpha,t_\alpha)/u =
(\alpha,\alpha)/u$, where $t_\alpha'$ is the element of $H$ corresponding to
$\alpha$ defined by ${\Cal B}'$. Modulo some scaling, we may assume that
$(\alpha,\alpha)\in\Bbb Q$ for one $\alpha\in R^\times$, multiplied by a nonzero scalar, e.g.
$u=(\alpha,\alpha)$. If $\beta\in R^\times$, we have
$\frac{2(\beta,\alpha)}{(\alpha,\alpha)}\in\Bbb Z$ so that $(\beta,\alpha)\in\Bbb Q$, and hence, since
$\frac{2(\alpha,\beta)}{(\beta,\beta)}\in\Bbb Z$, we get
$(\beta,\beta)\in\Bbb Q$ if $(\alpha,\beta)\neq 0$. Thus, by (A5), our form is $\Bbb Q$-valued on the
$\Bbb Q$-linear span of anisotropic roots. From now on {\it we assume that our form is scaled, if
$R^\times\neq\emptyset$, so that there is at least one $\alpha\in R^\times$ with
$(\alpha,\alpha)>0$, and that $(\beta,\beta) \in \Bbb Q$ for all $\beta \in R^\times$.} But then we can
immediately prove the following (shown in [A-P. Lem.I.2.3]).

\proclaim{Lemma 3.1} Let $\gamma\in R^\times$. Then $(\gamma,\gamma)>0$.
\endproclaim We give an elementary proof which is different from the one in [A-P].

\demo{Proof} Suppose not, i.e., $(\gamma,\gamma)<0$. Then, by (A5), there exists $\alpha,\beta\in
R^\times$ such that $(\alpha,\alpha)> 0$,
$(\beta,\beta)< 0$ and $(\alpha,\beta)\neq 0$. By Lemma 2.2, we have
$\alpha+\beta\notin R$ or $\alpha-\beta\notin R$. If $\alpha+\beta\notin R$, then the root string
$\beta-d\alpha,\ldots,\beta-\alpha,\beta$ with
$d=\frac{2(\beta,\alpha)}{(\alpha,\alpha)}>0$ gives $(\beta,\alpha)>0$. However, the root string
$\alpha-d'\beta,\ldots,\alpha-\beta,\alpha$ with
$d'=\frac{2(\alpha,\beta)}{(\beta,\beta)}>0$ gives $(\alpha,\beta)<0$, which is a contradiction. If
$\alpha-\beta=\alpha+(-\beta)\notin R$, we also get a contradiction since
$(-\beta,-\beta)< 0$ and $(\alpha,-\beta)\neq 0$.
\qed
\enddemo We later use the following. The proof is the same as in [A-P, Lem.I.2.6] using Lemma 2.1(5)
and Lemma 3.1:
\proclaim{Lemma 3.2} Let $\alpha \in R^\times$ and $\xi \in R$. Then
$-4 \leq \xi(h_\alpha) \leq 4$.
\endproclaim

Now, we shall prove a crucial property that the isotropic roots are in the radical of the form. This was
proved in [A-P, Prop.I.2.1] using two extra assumptions, namely,
$R$ is discrete in $H^*$ (assuming the base field is $\Bbb C$) and
$R^0$ is not isolated. We do not need such assumptions.  The first part of our proof has already been
established in in a few lines in the recent preprint [AKY] (in a different setup). It turns out that a
small modification of the proof in [A-P] was enough to exclude the extra assumptions.

We start with a lemma [A-P, Lem.I.1.30].  (This is well-known in terms of the representation theory for
the 3-dimensional Heisenberg Lie algebra. Namely, there is no nontrivial finite-dimensional
representation of the Lie algebra. [B])

\proclaim{Lemma 3.3} Let $\delta\in R^0$ and $\xi\in R$. Suppose that $(\xi,\delta)\neq 0$. Then
$\xi+n\delta\in R$ for infinitely many integers $n$.
\endproclaim

For the convenience of the reader, and also since [AKY] starts with different axioms, we repeat their
argument  to make sure that our axioms are enough. Also, note that it works for a null system.

\proclaim{Proposition 3.4} Let $(L,H,{\Cal B}) \in \Cal L$ be a LEALA. Let $\xi\in R$ and
$\delta\in R^0$. Then $(\xi,\delta)=0$. That is, $(R,R^0)=0$.
\endproclaim
\demo{Proof} First we show the statement for $\alpha:=\xi\in R^\times$. Suppose that
$(\alpha,\delta)\neq 0$. Then, by Lemma 3.3, we have $\alpha+n\delta\in R$ for infinitely many integers
$n$. But
$\alpha+n\delta$ can be isotropic for at most one $n$ (since $(\alpha+n\delta,\alpha+n\delta)=0$
$\Rightarrow$ $(\alpha,\alpha)+2n(\alpha,\delta)=0$
$\Rightarrow$ $n=-\frac{2(\alpha,\delta)}{(\alpha,\alpha)}$). Hence,  $\alpha+n\delta\in R^\times$ for
infinitely many integers $n$, and by Lemma 2.1(1),
$$
\frac{2(\delta,\alpha+n\delta)}{(\alpha+n\delta,\alpha+n\delta)}
=\frac{2(\delta,\alpha)}{(\alpha,\alpha)+2n(\alpha,\delta)}
\in\Bbb Z
$$ for infinitely many integers $n$, which is impossible. Therefore, $(\alpha,\delta)= 0$.

Suppose that $\delta':=\xi\in R^0$ and $(\delta',\delta)\neq 0$. Then,
$(\delta \pm \delta',\delta \pm \delta') =
\pm 2(\delta,\delta') \not= 0$. If $\delta+\delta' \in R$ (resp. $\delta-\delta' \in R$), then
$\delta+\delta' \in R^\times$ (resp. $\delta-\delta' \in R^\times$). In this case, we have $0 =
(\delta,\delta \pm
\delta') =
\pm (\delta,\delta') \not= 0$, which is a contradiction. If $\delta \pm \delta' \not\in R$, then
$[x_\delta,x_{\delta'}] = [y_\delta,x_{\delta'}] = 0$, and so $[t_\delta,x_{\delta'}] = 0$. (The
notations are from Heisenberg triplets, defined in \S 1.) This implies
$[t_\delta,x_{\delta'}] = (\delta,\delta')x_{\delta'}$, and hence $(\delta,\delta') = 0$, a
contradiction again.
\qed
\enddemo

Let $(L,H,{\Cal B}) \in \Cal L$ be a LEALA or a null system. Let $V$ be the $\Bbb Q$-subspace of
$H^\ast$ spanned by $R$. We will prove that the form on $V$ is positive semidefinite (Kac conjecture).
First as an immediate corollary of Proposition 3.4, we have:
\proclaim{Corollary 3.5} Let $(L,H,{\Cal B}) \in \Cal L$ be a null system. Then the form $(\cdot,\cdot)$ on $V$ is zero.

\endproclaim Let $(L,H,{\Cal B}) \in \Cal L$ be a LEALA. Put 
$V^0 = \{ v \in V \mid (v,w) = 0\ \text{\rm for\ all}\ w \in V \}$,  the radical of $(\cdot,\cdot)$, and
$\bar V = V/V^0$. Let
$\bar { } : V \longrightarrow \bar V$ be the canonical map. Then the form $(\cdot,\cdot)$ on $V$
induces a unique nondegenerate symmetric bilinear form $(\cdot,\cdot)$ on
$\bar V$ so that $(\bar v,\bar w) = (v,w)$ for all $v,w \in V$. Define $\bar R = \{ \bar \xi \mid \xi
\in R \}$ and
$\bar R^\times = \{ \bar \alpha \mid \alpha \in R^\times \}$. Then $\bar R = \bar R^\times \cup \{
\bar 0 \}$. Hence, $(\bar \alpha,\bar \alpha) > 0$ for all
$\bar \alpha \in \bar R$ with $\bar \alpha \not= \bar 0$, and $2(\bar \beta,\bar \alpha)/(\bar
\alpha,\bar \alpha)
\in {\Bbb Z}$.

We note that the following result is true in general, which is usually discussed only for the
finite-dimensional case  in standard textbooks.

\proclaim{Lemma 3.6} Let $X$ be a vector space over any field  with nondegenerate symmetric bilinear
form
$(\cdot, \cdot)$. Let $Y$ be a finite-dimensional subspace of $X$. Then there exists a
finite-dimensional subspace $\tilde Y$ of $X$ containing $Y$ such that $(\cdot,\cdot)\mid_{\tilde
Y\times\tilde Y}$ is nondegenerate.
\endproclaim
\demo{Proof} Let $\{ y_1, \ldots, y_r, u_1, \ldots, u_m\}$ be a basis of $Y$ so that
$\{ u_1, \ldots, u_m\}$ is a basis of the radical of $(\cdot,\cdot)\mid_{Y\times Y}$. If $m=0$ (i.e.,
the radical is zero), we can take $\tilde Y=Y$, and so we assume $m>0$. Since $(\cdot, \cdot)$ is
nondegenerate on $X$, there exists $x_1\in X$ such that $(u_1,x_1)=1$. Let $Y_1:=\langle y_1, \ldots,
y_r, u_1, x_1\rangle$. Let 
$y_i':=y_i-(y_i, x_1)u_1$ for $i=1,\ldots, r$. Then
$Y_1=\langle y_1', \ldots, y_r'\rangle\perp\langle u_1, x_1\rangle$, and so
$(\cdot,\cdot)\mid_{Y_1\times Y_1}$ is nondegenerate. Let $\tilde Y_1:=Y_1+ \langle u_2,\ldots,
u_m\rangle$. Let
$u_j':=u_j-(x_1, u_j)u_1$ for $j=2,\ldots, m$. Then
$Y\subset \tilde Y_1 =Y_1\perp\langle u_2',\ldots, u_m'\rangle$, and the radical
$\langle u_2',\ldots, u_m'\rangle$ on $\tilde Y_1$ has dimension $m-1$. Hence, by induction, there
exists a subspace $\tilde Y$ of $X$ containing $\tilde Y_1$ such that $(\cdot,\cdot)\mid_{\tilde
Y\times\tilde Y}$ is nondegenerate.
\qed
\enddemo

Using this lemma and a technique similar in [A-P, Lem.I.2.10], we can prove something new:

\proclaim{Proposition 3.7}
$(\bar R, \bar V)$ is a locally finite irreducible root system (in the sense of [LN]).  Also, $\bar R$ is finite if $\dim_F
H<\infty$ (i.e., if $L$ is an EALA), and hence $\bar V$ is finite-dimensional in this case, and so $(\bar R, \bar V)$ is a finite irreducible root system (in the sense of [Bo, Ch.VI]).
\endproclaim

\demo{Proof} We choose a basis $\{\bar\alpha_i\}_{i\in I}$ of $\bar V$. Let  ${\bar W}$ be a
finite-dimensional subspace of $\bar V$. Suppose that $\bar R\cap \bar W\neq 0$. Then, by Lemma 3.6,
there exists a finite-dimensional subspace $\Cal U$ of $\bar V$ containing
${\bar W}$ such that
$(\cdot,\cdot)\mid_{\Cal U\times\Cal U}$ is nondegenerate. Let $\{\bar w_1,\ldots,\bar w_r\}$ be a
basis of $\Cal U$. Since  $(\cdot,\cdot)\mid_{\Cal U\times\Cal U}$ is nondegenerate,  the map 
$f:\bar \beta \mapsto ((\bar w_1,\bar \beta^\vee),\ldots, (\bar w_r,\bar \beta^\vee))$ of $\bar R\cap
{\bar W}$ into
$\Bbb Q^r$ is injective, where
$\bar\beta^\vee = 2\bar \beta/(\bar \beta,\bar \beta)$. Since 
$\bar w_i=\sum_{j\in I_i}\ a_{ij}\bar\alpha_j$ for all $1\leq i\leq r$, some  finite subset $I_i$ of
$I$ and
$a_{ij}\in\Bbb Q$, we have, by Lemma 3.2,
$$ (\bar w_i,\bar \beta^\vee) =\sum_{j\in I_i}\ a_{ij}(\bar\alpha_j,\bar \beta^\vee)=
\sum_{j\in I_i}\ a_{ij}\alpha_j(h_{\beta})\in
\{\sum_{j\in I_i}\ m_ja_{ij}\mid -4\leq m_j\leq 4,\ m_j\in\Bbb Z\}.
$$ Hence the image of $f$ is finite and so is
$\bar R\cap {\bar W}$, i.e., $\bar R$ is locally finite. By Lemma 2.1(2), we have
$\sigma_{\bar \alpha}(\bar R) =\bar R$, and so $(\bar R, \bar V)$ is a locally finite root system
relative to the map ${}^\vee:\bar R\longrightarrow \bar V^*$ defined by the pairing
$\langle \bar v,\bar{\alpha}\rangle:=(\bar v,\bar{\alpha}^\vee)$ for all $\bar v\in \bar V$ and
$\bar{\alpha}\in\bar R$. The irreducibility follows from (A5).

For the second statement, let $V_F$ be the subspace of $H^*$ generated by $R$ over $F$. Put 
$V_F^0 = \{ v \in V_F \mid (v,w) = 0\ \text{for all}\ w \in V_F\}  = \{ v \in V_F \mid (v,\beta) = 0\
\text{for\ all}\ \beta \in R \}$,  the radical of $(\cdot,\cdot)$ on $V_F$. Then we see that $V^0 = V
\cap V_F^0$. Therefore, 
$\bar V_F:=V_F/V_F^0 \supset V/V^0=\bar V\supset\bar R$.  Thus if $H^*$ is finite-dimensional, then so
is $V_F$ or $\bar V_F$. We choose an $F$-basis $\{ \bar\alpha_1,\ldots,\bar\alpha_\ell \}$ of $\bar
V_F$ from $\bar R$. Since $(\cdot,\cdot)$ is nondegenerate on $\bar V_F$, the map
$\bar \beta \mapsto ((\bar \beta,\bar \alpha_1^\vee),\ldots, (\bar \beta,\bar \alpha_\ell^\vee))$ of
$\bar R$ into
$\{ -4,\ldots,0,\ldots,4 \}^\ell$ is injective (by Lemma 3.2). Hence $\bar R$ is finite.
\qed
\enddemo

One of the important notions for subsets of a locally finite root system is the fullness (see [LN] or
[NS]).

\definition{Definition 3.8} Let $\Delta$ be a locally finite root system and $M$ a subset of $\Delta$.
Then
$\Delta_M:=(\spa_\Bbb Q M)\cap \Delta$ is called the   {\it full subsystem generated by $M$}.

\enddefinition

The following lemma [LN, Cor.3.16] is useful.
\proclaim{Lemma 3.9} Let $M$ be a finite subset of a locally finite irreducible root system $\Delta$.
Then the full subsystem $\Delta_M$ is finite and irreducible.
\endproclaim

The Kac conjecture follows from [LN, Thm 4.2] or from an argument analogous to the one given in [A-P]:
Let $0\neq
\bar v\in\bar V$. Then $\bar v=\sum_{\bar\alpha_i\in M}\ a_i\bar\alpha_i$ for $a_i\in\Bbb Q$ and some
finite subset
$M$ of $\bar R$. Then, by Lemma 3.9, the full subsystem $\bar R_M$ is finite and irreducible so that 
$\bar W:=\spa_\Bbb Q M$ contains $\bar v$. Then, $(\cdot,\cdot)\mid_{\bar W\times \bar W}$ is
nondegenerate since the Cartan matrix of $\bar R_M$ is nonsingular. Then, apply for $(\bar R_M,\bar W)$
instead of $(\bar R,\bar V)$ in [A-P, Thm I.2.14] to get
$(\bar v,\bar v)>0$ (since
$(\bar \alpha,\bar \alpha) > 0$ for all
$\bar \alpha \in \bar R_M$).  Thus:
\proclaim{Theorem 3.10} Let $(L,H,{\Cal B}) \in \Cal L$ be a LEALA. Then the form $(\cdot,\cdot)$ on
$V$ is positive semidefinite.
\endproclaim

\definition{Definition 3.11} The dimension of the radical for $V$ is called the {\it nullity} for a
LEALA. If the additive subgroup of $V$ generated by $R^0$ is free, we call the rank the {\it null rank}
for a LEALA.

\enddefinition

The null rank and the nullity coincide for the definition in [A-P] or [N2]. In our definition, null
rank $n$ implies nullity $n$, but nullity $n$ does not imply null rank $n$. We will give an example of
an EALA of nullity 1, which does not have any null rank in \S 5.

\head
\S 4 Root systems
\endhead

One can naturally consider a generalization of extended affine root systems defined by Saito [S], and
the anisotropic root systems of our LEALAs  (in $V\otimes_\Bbb Q\Bbb R$)  are examples of the
generalized root systems, defined in the following.

\definition{Definition 4.1} Let $V$ be a vector space over $\Bbb R$ with (nonzero) positive
semidefinite form
$(\cdot,\cdot)$. A subset $\frak R$ of $V$ is called a {\it locally extended affine root system} or a
{\it LEARS} for short if
\roster
\item $(\alpha,\alpha)\neq 0$  for all $\alpha\in\frak R$, and
$\frak R$ generates $V$;

\item
$\frac{2(\alpha,\beta)}{(\alpha,\alpha)}\in\Bbb Z$ for all 
$\alpha, \beta\in\frak R$;

\item
$\sigma_\alpha(\frak R)\subset \frak R$ for all $\alpha\in \frak R$, where
$\sigma_\alpha(\beta)=\beta-\frac{2(\alpha,\beta)}{(\alpha,\alpha)}\alpha$;

\item $\frak R=\frak R_1\cup \frak R_2$ and $(\frak R_1,\frak R_2)=0$ imply $\frak R_1=\emptyset$ or
$\frak R_2=\emptyset$. (Or $\frak R$ is irreducible.)

\endroster

\enddefinition

Let $(V,\frak R)$ be a LEARS. Let $(\bar V,\bar \frak R)$ be the canonical image onto $V$ modulo the
radical as in
\S 3. Then $\bar V$ admits the positive definite form, and so as in [A-P, Lem.II.2.8], we get the
following.
\proclaim{Proposition 4.2} Let $(V,\frak R)$ be a LEARS. Then $(\bar V,\bar \frak R)$ is a locally
finite irreducible root system.
\endproclaim
\demo{Proof} By Schwarz' inequality, we have
$$\left |\matrix 
\frac{2(\bar\beta,\bar\alpha)}{(\bar\alpha,\bar\alpha)}  & 
\!\!\!\!\frac{2(\bar\alpha,\bar\beta)}{(\bar\beta,\bar\beta)}
\endmatrix
\right |
\leq 4
$$ for all $\bar\alpha,\bar\beta\in\bar \frak R$. Thus, for $\alpha,\beta \in \frak R$, we get $-4 \leq
\frac{2(\beta,\alpha)}{(\alpha,\alpha)} \leq 4$, and repeat the argument in Proposition 3.7.
\qed
\enddemo

\example{Remark 4.3} If we define the set
$S_{\bar\alpha}=\{\delta\in V\mid \ol{\alpha+\delta}\in\bar\frak R\}$ for $\bar\alpha\in\bar\frak R$,
then
$$\frak R=\bigcup_{\bar\alpha\in\bar\frak R}\ (\alpha+S_{\bar\alpha}).$$ Let ${\Lambda}$ be the
additive subgroup of $V$ generated by
$\cup_{\bar\alpha\in\bar\frak R}\ S_{\bar\alpha}$. Then the family of subsets $S_{\bar\alpha}$ of
${\Lambda}$, say
$\{S_{\bar\alpha}\}_{\bar\alpha\in\bar\frak R}$, is a natural generalization of root systems extended by
${\Lambda}$, introduced in [Y2]. (A finite irreducible root system was taken as $\bar\frak R$ in [Y2],
and here a locally finite irreducible root system is taken as $\bar\frak R$.)

\endexample

\head
\S 5 Examples of new EALAs
\endhead

First we construct an analogue of loop algebras. Let ${\Lambda}=({\Lambda},+,0)$ be an abelian group.
Let
$$F[{\Lambda}] =\bigoplus_{{\lambda}\in {\Lambda}}\ F\ol {\lambda}$$ be the group algebra of
${\Lambda}$ over $F$. We define a bilinear form
$\varepsilon$ on $F[{\Lambda}]$ by
$$
\varepsilon(\ol {\lambda}, \ol {\mu}):=
\cases 1&\text{if ${\lambda}+{\mu}=0$}\\ 0&\text{otherwise}
\endcases
$$ for
${\lambda},{\mu}\in {\Lambda}$. Then $\varepsilon$ is a nondegenerate symmetric invariant form. Let
$\frak g$ be a finite-dimensional split simple Lie algebra of type $\Delta$ with a split Cartan
subalgebra $\frak h$. Let
$$M:=\frak g\otimes_F F[{\Lambda}]$$ be the Lie algebra
 with the bracket defined by
$$[x\otimes\ol {\lambda}, y\otimes\ol {\mu}]=[x,y]\otimes \ol{{\lambda}+{\mu}}$$ for $x,y\in\frak g$
and ${\lambda},{\mu}\in {\Lambda}$. Note that if ${\Lambda}=\Bbb Z$, then $M$ is a loop algebra. Let
$$(\cdot,\cdot):=\kappa\otimes\varepsilon,$$ where $\kappa$ is the Killing form of $\frak g$. Then
$(\cdot,\cdot)$ is a nondegenerate symmetric invariant bilinear form.

We assume that there exists a nonzero additive homomorphism $\vvp$ of ${\Lambda}$ into $F$. Let
$d_{\vvp}$ be a degree derivation of $F[{\Lambda}]$ determined by $\vvp$, i.e.,
$d_{\vvp}(\ol {\lambda})=\vvp({\lambda})\ol {\lambda}$ for
${\lambda}\in {\Lambda}$. We define a Lie algebra $L$ by
$$ L=M\oplus Fc\oplus Fd_\vvp,
$$ where $c$ is a nonzero central element with multiplication as follows:
$$
\align [d_\vvp,x\otimes \ol {\lambda}] &=x\otimes d_\vvp(\ol {\lambda}) =-[x\otimes \ol {\lambda},
d_\vvp]\quad\text{and}
\tag5.1\\ [x\otimes \ol {\lambda},y\otimes \ol {\mu}]&=[x,y]\otimes 
\ol{{\lambda}+{\mu}}+ (x\otimes d_\vvp(\ol {\lambda}),y\otimes \ol {\mu})c
\endalign
$$ for all $x\otimes \ol {\lambda},y\otimes \ol {\mu}\in M$. (Indeed, this is a Lie algebra since
$d_\vvp$ is a skew derivation relative to $\varepsilon$.) Also, one can extend the form $(\cdot,\cdot)$
to a form
${\Cal B}(\cdot,\cdot)$ on $L$ by
$$ {\Cal B}(c,c)={\Cal B}(d_{\vvp},d_{\vvp})={\Cal B}(c,M)={\Cal B}(d_{\vvp},M)=0
\quad\text{and}\quad {\Cal B}(c,d_{\vvp})=1.
$$ Then
${\Cal B}(\cdot,\cdot)$ is a nondegenerate symmetric invariant bilinear form. Let
$$H=\frak h\oplus Fc\oplus Fd_{\vvp}.$$ Let $\delta_{\lambda}$ for ${\lambda}\in {\Lambda}$ be  the
element of $H^*$ such that
$\delta_{\lambda}(d_\vvp)=\vvp({\lambda})$ and $\delta_{\lambda}(\frak h)=\delta_{\lambda}(c)=0$. Then
$$[d_\vvp,h\otimes \ol {\lambda}]=d_\vvp({\lambda})h\otimes \ol {\lambda}  =\vvp({\lambda})h\otimes \ol
{\lambda} =\delta_{\lambda}(d_\vvp)h\otimes \ol {\lambda}$$ and
$[\frak h\oplus Fc,h\otimes\ol {\lambda}]=0$ for all $h\in\frak h$.  So $\delta_{\lambda}$ is a root
relative to $H$. Note that the centralizer of $H$ is contained in
$(\frak h\otimes F[{\Lambda}])\oplus Fc\oplus Fd_{\vvp}$. Thus $H$ is self-centralizing if and only if
$\vvp$ is injective. We now assume that $\vvp$ is injective so that $(L,H,{\Cal B})$ is an admissible
triple. Note that
$$\text{${\Lambda}$ has to be torsion-free.}
$$ Then one can check that if we denote the root system of $\frak g$ by $\Delta$, then
$$R^0=\{\delta_{\lambda}\mid {\lambda}\in {\Lambda}\}
\quad\text{and}\quad R^\times=\{\alpha+\delta_{\lambda}\mid \alpha\in\Delta,\ {\lambda}\in {\Lambda}\},
$$ and $L$ is an EALA. In particular, if we take ${\Lambda}$ to be any additive subgroup of $F$ and
$\vvp$ to be the inclusion, then
$R^0\cong {\Lambda}$ and the nullity of $L$ is equal to the dimension of the $\Bbb Q$-span of
${\Lambda}$ over $\Bbb Q$. For example, if ${\Lambda}=\Bbb Q$, then
$L$ is an EALA of nullity $1$ and $R^0=\langle R^0\rangle$ is not a free abelian group (and so
$\langle R\rangle$ is not a lattice). So  $L$ does not have null rank.

To give an example of null rank $\infty$ (so nullity $\infty$),  let ${\frak I}$ be an index set of an
arbitrary cardinality. We assume that $F$ contains a linearly independent subset $S$ over $\Bbb Q$ with
$|S|=|{\frak I}|$. Let ${\Lambda}=\Bbb Z^{\oplus {\frak I}}$. Then there is a group isomorphism  from
${\Lambda}$ into the additive subgroup of $F$ generated by $S$, say $\vvp$. Thus our construction above
gives an EALA with
$R^0=\Bbb Z^{\oplus {\frak I}}$ and the null rank is $|{\frak I}|$. 

We note that this example is a natural generalization of untwisted affine Kac-Moody Lie algebras from
the loop algebra $\frak g\otimes F[t^{\pm 1}]$ to an infinite-loop algebra  $$\frak g\otimes F[t_i^{\pm 1}]_{i\in\Bbb N}.$$   where $F[t_i^{\pm 1}]_{i\in\Bbb N}$ is the ring of Laurent polynomials in infinitely many variables. To see this more clearly, we give a slightly different description.

Let $F=\Bbb C$ for convenience. (One can take $F$ to be $\Bbb R$ or any field with $[F:\Bbb
Q]=\infty$.) In fact,
$$\Bbb C[t_i^{\pm 1}]_{i\in\Bbb N} =\bigoplus_{\bs\alpha\in\Bbb Z^{\oplus\Bbb N}}\ \Bbb
Ct_{\bs\alpha},$$ where $t_{\bs\alpha}=\prod_{i\in \Bbb N}\ t_i^{\alpha_i}$ for
$\bs\alpha=(\alpha_i)\in\Bbb Z^{\oplus \Bbb N}$ (only finite number of $\alpha_i$ are not zero),
 is the group algebra $\Bbb C[\Bbb Z^{\oplus \Bbb N}]$ by $\ol{\bs\alpha}=t_{\bs\alpha}$. Let
$\{s_i\}_{i\in\Bbb N}$ be a linearly independent subset of $\Bbb C$ over $\Bbb Q$, and let
$$d:=\sum_{i\in\Bbb N}\ s_it_i{\partial\over\partial t_i}
\quad\text{(formal infinite sum)}$$  be the derivation of $\Bbb C[t_i^{\pm 1}]_{i\in\Bbb N}$, which can
be interpreted as our $d_\vvp$ above. Thus the second part of (5.1) can be rephrased as follows:
$$ [x\otimes t_{\bs\alpha},y\otimes t_{\bs\beta}] =[x,y]\otimes t_{\bs\alpha+\bs\beta}+
\kappa(x,y)\sum_{i\in\Bbb N}\ s_i\text{Res}_i({\partial t_{\bs\alpha}\over\partial t_i} t_{\bs\beta})c
$$ for all $x\otimes t_{\bs\alpha},y\otimes t_{\bs\beta}\in M$, where
$\text{Res}_i$ gives the coefficient of
$t_i^{-1}$.

\example{Remark 5.2} (1) In the setting above, let
$|I|=2$,
$s_1=1$ and $s_2=\sqrt 2$ (or any irrational number). Then
$R^0=\Bbb Z^2$ and so the null rank is 2, but
$R^0$ is not discrete in the $\Bbb R$-span of $R^0$. So this is not an EALA in the sense of [A-P] (see
also [G]). But in our sense this is just an EALA of null rank 2 over $\Bbb C$.

(2) The algebra
$\Bbb C[t_i^{\pm 1}]_{i\in\Bbb N}$ can be generalized to a quantum torus of infinitely many variables
$\Bbb C_{\bs q}[t_i^{\pm 1}]_{i\in\Bbb N}$, and from 
$sl_n(\Bbb C_{\bs q}[t_i^{\pm 1}]_{i\in\Bbb N})$, we get an EALA of null rank $\infty$ by the same
construction. More generally, from Jordan or structurable ${\Lambda}$-tori ([Y1], [AY] or [Y3]), where
${\Lambda}$ is a torsion-free abelian group, one can construct various new EALAs. 

\endexample

\head
\S 6 Examples of LEALAs
\endhead

Let $\frak g_\infty$ be one of the so-called infinite rank affine algebras of type $\text A_{+\infty}$,
$\text A_{\infty}$, $\text B_\infty$, $\text C_\infty$ and $\text D_\infty$, and identify it with a
subalgebra of $gl_\infty(F)$, the Lie algebra of
all matrices $(a_{ij})_{i,j\in\Bbb Z}$ ($a_{ij}\in F$) such that the number of nonzero $a_{ij}$
is finite, with the usual bracket, as in [K, \S 7.11]. 
Then it is easy to show that $\frak g_\infty$  with
the trace form and taking our Cartan to be the subalgebra consisting of diagonal matrices is a LEALA of
nullity 0. Also, as in \S 5, one can construct a LEALA 
$$(\frak g_\infty\otimes_\Bbb C \Bbb C[t_1^{\pm 1},\ldots , t_n^{\pm 1}])
\oplus \Bbb C c\oplus \Bbb C d_\vvp
\quad\text{(null rank $n$) or}$$
$$(\frak g_\infty\otimes_\Bbb C \Bbb C[t_i^{\pm 1}]_{i\in\Bbb N})
\oplus \Bbb C c\oplus \Bbb C d_\vvp
\quad\text{(null rank $\infty$),}$$ taking the trace form instead of the Killing form for the
multiplication. Moreover, as in Remark 5.2(2),
$$sl_\infty(\Bbb C_{\bs q}[t_i^{\pm 1}]_{i\in\Bbb N})
\oplus \Bbb Cc\oplus \Bbb Cd_\vvp$$ is a LEALA of null rank $\infty$. ($sl_\infty(A)$ for any
associative algebra
$A$ is defined in [N1].)

Locally finite split simple Lie algebras were classified by Neeb and Stumme in [Stu] and [NS]. They
showed that locally finite split simple Lie algebras over a field of characteristic 0 are isomorphic to
exactly one of infinite rank affine Lie algebras of type $\text A_{\infty}$, $\text B_\infty$ and
$\text C_\infty$. (They showed that $\text A_{+\infty}\cong \text A_{\infty}$ and  $\text
B_{\infty}\cong \text D_{\infty}$.) To classify LEALAs of nullity 0 in \S 8, we need more information
about these algebras.  Thus we precisely define them here.
(The size of matrices is not only $\aleph_0$ but any cardinality.) 
\definition{Definition 6.1} Let ${\frak I}$ be any index set. Then the Lie algebra of type $\text
X_{\frak I}$ is defined as a subalgebra of the matrix algebra $gl_{\frak I}(F)$,
$gl_{2{\frak I}+1}(F)$ or $gl_{2{\frak I}}(F)$ consisting of matrices having only finite nonzero
entries:

Type $\text A_{\frak I}$;
$sl_{\frak I}(F)=\{x\in gl_{\frak I}(F)\mid \tr (x)=0\}$;

Type $\text B_{\frak I}$;
$o_{2{\frak I}+1}(F)=\{x\in gl_{2{\frak I}+1}(F)\mid sx=-x^ts\}$;

Type $\text C_{\frak I}$;
$sp_{2{\frak I}}(F)=\{x\in gl_{2{\frak I}}(F)\mid sx=-x^ts\}$;

Type $\text D_{\frak I}$;
$o_{2{\frak I}}(F)=\{x\in gl_{2{\frak I}}(F)\mid sx=-x^ts\}$;

\noindent where 
$$s=
\left (\matrix \format\r&\quad\r &\quad\r    \\
   0   &I_{\frak I}    & 0   \\
  I_{\frak I}   & 0    & 0  \\
   0   & 0   &   1   
\endmatrix
\right ),\
\left (\matrix \format\r&\quad\r &\quad\r &\quad\r   \\
   0   &I_{\frak I}       \\
  -I_{\frak I}   & 0   \\
\endmatrix
\right )
\quad\text{or}\quad
\left (\matrix \format\r&\quad\r &\quad\r &\quad\r   \\
   0   &I_{\frak I}       \\
  I_{\frak I}   & 0   \\
\endmatrix
\right ),
$$ respectively for type $\text B_{\frak I}$, $\text C_{\frak I}$ or $\text D_{\frak I}$, and $I_{\frak
I}$ is the identity matrix of size ${\frak I}$ ($x^t$ = transpose of $x$).

\enddefinition

We give interesting examples of LEALAs of nullity 0.

\example{Example 6.2} (1) Let $d=\sum_{i\in\Bbb N}\ ie_{ii}$, where the $e_{ij}$ are matrix units in
$M_{\Bbb N}(F)$. Define a bilinear form ${\Cal B}$ on the matrix Lie algebra $L=sl_{\Bbb N}(F)\oplus
Fd$ by
${\Cal B}(x,y):=\tr(xy)$ for $x\in sl_{\Bbb N}(F)$ and $y\in L$, and ${\Cal B}(d,d):=0$. Let $\frak h$
be the subalgebra of $sl_{\Bbb N}(F)$ consisting of diagonal matrices, and let $H:=\frak h\oplus Fd$.
Then
$(L,H,{\Cal B})$ is a centreless LEALA of nullity 0.

(2) Let $d=\sum_{i\in\Bbb N}\ (e_{ii}-e_{\Bbb N+i,\Bbb N+i})$. Define a bilinear form ${\Cal B}$ on the
natural Lie algebra $L=o_{2\Bbb N+1}(F)\oplus Fd$ by ${\Cal B}(x,y):=\tr(xy)$ for $x\in o_{2\Bbb
N+1}(F)$ and
$y\in L$, and ${\Cal B}(d,d):=0$. Let $\frak h$ be the subalgebra of $o_{2\Bbb N+1}(F)$ consisting of
diagonal matrices, and let $H:=\frak h\oplus Fd$. Then $(L,H,{\Cal B})$ is a centreless LEALA of
nullity 0.

Note that one can construct similar kinds of Lie algebras of type $\text C_{\frak I}$ or $\text
D_{\frak I}$. Also, for any case, you can choose any scalar for ${\Cal B}(d,d)$ (not necessarily 0).

\endexample

Finally, we will show in Proposition 8.3 that the {\it core} of a LEALA (which will be defined in 8.1)
is a direct limit of Lie tori (which are defined in [N3] or [Y2]). Thus, to classify LEALAs, one may
need some ideas from the theory of locally finite Lie algebras (e.g. [BB], [GN] or [St]).

\head
\S 7 Examples of null systems
\endhead

A typical example of null systems is a Heisenberg Lie algebra with derivations added. More precisely:

\example{Example 7.1} Let ${\Lambda}=({\Lambda}, +, 0)$ be an abelian group. Let $S$ be a subset of
${\Lambda}$ satisfying
$$0\in S,
\quad\text{and}\quad
\delta\in S \Rightarrow -\delta\in S.
$$ Let $N=\oplus_{\delta\in S}\ N_\delta$ be a graded vector space over $F$ with a symmetric bilinear
form ${\Cal B}$ satisfying
\roster
\item"(N1)" ${\Cal B}\mid_{N_\delta\times N_{-\delta}}$ is nondegenerate for each $\delta\in S$;
\item"(N2)" 
$N_0=Fc+ Fd$ for some $c,d\in N_0$ with ${\Cal B}(c,c)={\Cal B}(d,d)=0$ and ${\Cal B}(c,d)=1$ (or
$N_0$ is a hyperbolic plane);
\item"(N3)" ${\Cal B}\mid_{N_\delta\times N_{\delta'}}=0$ unless $\delta'=-\delta$.
\endroster

We assume that there exists an injective additive homomorphism $\vvp$ from ${\Lambda}$ into $F$. So
${\Lambda}$ has to be torsion-free as in \S 7. Now we define the Lie bracket on $N=(N,S,{\frak
b},\vvp)$ as follows. For any
$0\neq\delta,\delta'\in S$, $x\in N_\delta$ and $y\in N$, we define
$$ [x,y]={\Cal B}(\vvp(\delta)x,y)c,
\quad [c,N]=0=[N,c]
\quad\text{and}\quad [d,x]=\vvp(\delta)x=-[x,d].
$$ Then one can check that $N$ is in fact a Lie algebra, and $(N,N_0,{\Cal B})$ is a null system. Note
that
$Fc$ is the centre of $N$, and $N=[N,N]\oplus Fd$. As in \S 5, the nullity of $N$ is the $\Bbb
Q$-dimension of the $\Bbb Q$-span of $\vvp ({\Lambda})$ in $F$, and the null rank of $N$ is  the rank
of ${\Lambda}$ if
${\Lambda}$ is free.

Note that this $N$ is indecomposable unless $S=0$. In fact, suppose that $N$ is decomposable. Then
$N=(L_1,H_1,{\Cal B}_1)\oplus (L_2,H_2,{\Cal B}_2)$, and
$H_1=F(c+ud)$ and $H_2=F(c-ud)$ for some nonzero $u\in F$ since $H_1\perp H_2$. Let $0\neq x\in
N_\delta$ for any
$0\neq\delta\in S$. Then $[c+ud,x]=u\vvp(\delta)x\neq 0$, and so $x\notin L_2$. Similarly, $x\notin
L_1$. But we have
$N_\delta=(N_\delta\cap L_1)\oplus(N_\delta\cap L_2)$, which implies
$N_\delta=0$, a contradiction.
\endexample

\example{Remark 7.2} (1) If ${\Lambda}=\Bbb Z$, $S=\{0,\pm 1\}$,  and $\vvp=\id$, then $[N,N]$ is a
Heisenberg Lie algebra (see [MP, \S 1.5]). If we also assume that $\dim_FN_1=\dim_FN_{-1}=1$, then $N$
is 4-dimensional and $[N,N]$ is the 3-dimensional Heisenberg Lie algebra. Moreover, if one takes an
$n$-th copy of this 4-dimensional $N$, i.e., $N^{(n)}:=N\oplus \cdots \oplus N$, then
$[N^{(n)},N^{(n)}]$ is usually called a Heisenberg Lie algebra of order $n$ (see [K, \S 2.9]).

(2)  Let $L=(\frak g\otimes F[t^{\pm 1}])\oplus Fc\oplus Fd$ be an (untwisted) affine Kac-Moody Lie
algebra, which is a special case in the previous section, and let $N_m:=\frak h\otimes t^m$ for $m\neq
0$ and $N_0:=Fc\oplus Fd$. Then the $\Bbb Z$-graded subalgebra $N$ of $L$ is a null system ($S=\Bbb Z$)
and $[N,N]$ is a $\Bbb Z$-graded Heisenberg Lie algebra.

(3) If $S=0$, then $N=N_0$ is an abelian Lie algebra, which is a decomposable null system of nullity 0.
Conversely, any abelian Lie algebra can be considered as a null system of nullity 0.
\endexample

The next example is a generalization of the standard `null part' of an (original) EALA of  maximal type
for the case of commutative associative coordinates.
\example{Example 7.3} Consider $F[t_i^{\pm 1}]_{i\in\Bbb N} =\oplus_{\bs\alpha\in\Bbb Z^{\oplus\Bbb
N}}\ Ft_{\bs\alpha}$  as in \S 5. Let
$d_i$ be the natural projection of $\Bbb Z^{\oplus \Bbb N}$ onto $\Bbb Z$ for $i\in\Bbb N$, i.e.,
$\bs\alpha=(\alpha_i)\mapsto \alpha_i$, and let
$$\Cal D:=\bigoplus_{i\in\Bbb N}\ Fd_i.$$ Let
$$W:=\bigoplus_{\bs\alpha\in\Bbb Z^{\oplus \Bbb N}}\ t_{\bs\alpha}\Cal D
\subset\der (F[t_i^{\pm 1}]_{i\in\Bbb N})$$ be a generalized Witt algebra over $F$ so that
$$(t_{\bs\alpha}d)(t_{\bs\beta})=d(\bs\beta)t_{\bs\alpha+\bs\beta}$$ for $d\in \Cal D$ and
$\bs\alpha,\bs\beta \in \Bbb Z^{\oplus\Bbb N}$. The Lie bracket satisfies
$$ [t_{\bs\alpha}d, t_{\bs\beta}d'] = t_{\bs\alpha+\bs\beta}
 (d(\bs\beta)d' - d'(\bs\alpha)d),
\tag1$$ and so $W$ is a $\Bbb Z^{\oplus \Bbb N}$-graded Lie algebra. Let 
$W_{\bs\alpha}:=t_{\bs\alpha}\Cal D$ for
$\bs\alpha\in \Bbb Z^{\oplus\Bbb N}$. Then $W_\bk 0$ is a self-centralizing  abelian ad-diagonalizable
subalgebra. Let $\{\delta_i\}_{i\in\Bbb N}$ be the dual set of $\{d_i\}_{i\in\Bbb N}$. Identifying 
$\Bbb Z^{\oplus \Bbb N}\subset \Cal D^*$ by
$\bs\alpha=(\alpha_i)\mapsto \sum\ \alpha_i\delta_i$ (finite sum), the set of roots for $(W,W_\bk 0)$
is $\Bbb Z^{\oplus \Bbb N}$, i.e.,
$$ [d, t_{\bs\beta}d'] = t_{\bs\beta}d(\bs\beta)d' =\bs\beta(d)t_{\bs\beta}d'.
\tag2$$ Let 
$$\Cal C:=\bigoplus_{i\in\Bbb N}\ F\delta_i\subset \Cal D^*,$$ and let 
$$Y=\bigoplus_{\bs\alpha\in\Bbb Z^{\oplus \Bbb N}}\ t_{\bs\alpha}\Cal C$$ be a $\Bbb Z^{\oplus \Bbb
N}$-graded vector space. For $c\in \Cal C$, we define
$$ (t_{\bs\alpha}c)(t_{\bs\beta}d)=\delta_{\bs\alpha+\bs\beta,\bk 0}c(d).
\tag3$$ Then, $Y\subset W^*$. Consider $W^*$ as a natural $W$-module. So we have
$$\align ((t_{\bs\alpha}d).(t_{\bs\beta}c))(t_{\bs\gamma}d') &=(t_{\bs\beta}c)(-[t_{\bs\alpha}d,
t_{\bs\gamma}d'])\\ &=(t_{\bs\beta}c)(t_{\bs\alpha+\bs\gamma} (d'(\bs\alpha)d - d(\bs\gamma)d'))\\
&=\delta_{\bs\alpha+\bs\beta+\bs\gamma,\bk 0} (d'(\bs\alpha)c(d) - d(\bs\gamma)c(d'))\\
&=(t_{\bs\alpha+\bs\beta}(c(d) c_{\bs\alpha}+ d(\bs\alpha+\bs\beta)c))(t_{\bs\gamma}d'),
\endalign
$$ where
$c_{\bs\alpha}\in \Cal C$ is defined as
$c_{\bs\alpha}(d')=d'(\bs\alpha)$. (We could simply write $\bs\alpha$ instead of
$c_{\bs\alpha}$ as in (2).) Hence,
$$[t_{\bs\alpha}d,t_{\bs\beta}c] =t_{\bs\alpha+\bs\beta}(c(d) c_{\bs\alpha}+ d(\bs\alpha+\bs\beta)c).
\tag4$$ Letting
$Y_{\bs\alpha}:=\ t_{\bs\alpha}\Cal C$,  we get a graded $W$-module
$Y=\oplus_{\bs\alpha\in\Bbb Z^{\oplus \Bbb N}}\ Y_{\bs\alpha}$.

Now, for $\bs\alpha\in\Bbb Z^{\oplus \Bbb N}$, let 
$$W_{\bs\alpha}':=\{t_{\bs\alpha}d\in W_{\bs\alpha}\mid d(\bs\alpha)=0\}
\quad\text{and}\quad W'=\bigoplus_{\bs\alpha\in\Bbb Z^{\oplus \Bbb N}}\ W_{\bs\alpha}'.
$$ Then, for $t_{\bs\alpha}d,t_{\bs\beta}d'\in W'$, since $d(\bs\alpha)=d'(\bs\beta)=0$, we have (see
(1))
$$d(\bs\beta)d'(\bs\alpha+\bs\beta)  -d'(\bs\alpha)d(\bs\alpha+\bs\beta)
=d(\bs\beta)d'(\bs\alpha)-d'(\bs\alpha)d(\bs\beta)=0,$$ and hence
$[t_{\bs\alpha}d,t_{\bs\beta}d']\subset W_{\bs\alpha+\bs\beta}'$. So $W'$ is a graded subalgebra of
$W$. Let
$$ Y_{\bs\alpha}':=\{(t_{\bs\alpha}c)|_{W'}
\mid c\in \Cal C\}.
$$ Then $Y'=\oplus_{\bs\alpha\in\Bbb Z^{\oplus \Bbb N}}\ Y_{\bs\alpha}'$ is also a $W'$-graded module
with
$$[t_{\bs\alpha}d,t_{\bs\beta}c] =t_{\bs\alpha+\bs\beta}(c(d) c_{\bs\alpha}+ d(\bs\beta)c),
\tag4$'$ $$ instead of (4). Note that
$W_{\bk 0}'=W_{\bk 0}'$ and
$Y_{\bk 0}'=Y_{\bk 0}$. Let
$$ N:=W'\oplus Y'=\bigoplus_{\bs\alpha\in\Bbb Z^{\oplus \Bbb N}}\  (W_{\bs\alpha}'\oplus
Y_{\bs\alpha}').
$$ Consider $Y$ as an abelian Lie algebra. Then, by the action (4$'$),
$N=\oplus_{\bs\alpha\in\Bbb Z^{\oplus \Bbb N}}\  N_{\bs\alpha}$, where
$N_{\bs\alpha}=W_{\bs\alpha}'\oplus Y_{\bs\alpha}'$,  becomes a $\Bbb Z^{\oplus \Bbb N}$-graded Lie
algebra. (In other words, it is the trivial extension of $W'$ by the abelian Lie algebra $Y'$.) Also, 
$N_{\bk 0}=W_{\bk 0}\oplus Y_{\bk 0}$ is a self-centralizing abelian ad-diagonalizable subalgebra. For
we have
$$[c, t_{\bs\alpha}d](t_{\bs\beta}d') =-t_{\bs\alpha}(c(d) c_{\bs\alpha})(t_{\bs\beta}d')
=-c(d)d'(\bs\alpha)\delta_{\bs\alpha+\bs\beta,\bk 0} =c(d)d'(\bs\beta)\delta_{\bs\alpha+\bs\beta,\bk
0}=0,$$ and hence
$[W_{\bs\alpha}',Y_{\bk 0}]=0$. So 
$Y_{\bk 0}$ is the centre of $N$. Since
$[d,t_{\bs\beta}c] =t_{\bs\beta}d(\bs\beta)c$ (see (4$'$)), the set of roots is $\Bbb Z^{\oplus \Bbb
N}$ as in the case of $W$. Note that
$$[N_{\bs\alpha},N_{-\bs\alpha}]=Fc_{\bs\alpha}
\tag5$$ for all $\bs\alpha\in\Bbb Z^{\oplus \Bbb N}$. In fact, for $t_{\bs\alpha}d,t_{-\bs\alpha}d'\in
W'$, since
$d(\bs\alpha)=d'(\bs\alpha)=0$, we have 
$[t_{\bs\alpha}d,t_{-\bs\alpha}d']=0$ (since $d(-\bs\alpha)d' -d'(\bs\alpha)d=0$). Hence,
$[W_{\bs\alpha}',W_{-\bs\alpha}']=0$ for all $\bs\alpha\neq\bk 0$, and also,
$[t_{\bs\alpha}d,t_{-\bs\alpha}c] =c(d) c_{\bs\alpha}$ (by (4$'$)). Thus (5) holds.

We define a symmetric bilinear form ${\Cal B}$ as follows:
$$ {\Cal B}(W',W')={\Cal B}(Y',Y')=0
\quad\text{and}\quad {\Cal B}(t_{\bs\alpha}c,t_{\bs\beta}d) =(t_{\bs\alpha}c)(t_{\bs\beta}d)
=\delta_{\bs\alpha+\bs\beta,\bk 0}c(d).
$$ Then one can easily show that ${\Cal B}$ is invariant. To show the  nondegeneracy of ${\Cal B}$, let 
$\bs\alpha^\perp:=\{d\in \Cal D\mid d(\bs\alpha)=0\}$ and pick some $d_{\bs\alpha}\in \Cal D$ so that
$d_{\bs\alpha}(\bs\alpha)\neq 0$. Then $\Cal D=Fd_{\bs\alpha}+\bs\alpha^\perp$. Now,
$W_{\bs\alpha}'=\{t_{\bs\alpha}d\in W_{\bs\alpha}\mid d\in\bs\alpha^\perp\}$ and
$$ Y_{-\bs\alpha}'=\{(t_{-\bs\alpha}c)|_{W'}
\mid c\in \Cal C\} =\{(t_{-\bs\alpha}c)|_{W_{\bs\alpha}'}
\mid c\in \Cal C\} =\{t_{-\bs\alpha}c|_{\bs\alpha^\perp}
\mid c\in \Cal C\}.
$$ Identifying $c|_{\bs\alpha^\perp}$ with
$\tilde c\in \Cal C$, where
$\tilde c(d_{\bs\alpha})=0$ and $\tilde c|_{\bs\alpha^\perp}=c|_{\bs\alpha^\perp}$,
$$
\bs\alpha^\perp\times  \{c|_{\bs\alpha^\perp}
\mid c\in \Cal C\}
\longrightarrow F
$$ is nondegenerate (since the pairing
$\Cal D\times  \Cal C\longrightarrow F$ is nondegenerate). So $W_{\bs\alpha}'$ and $Y_{-\bs\alpha}'$
give a nondegenerate pair, and ${\Cal B}$ is nondegenerate. (We also see that
$W_{\bs\alpha}'\cong Y_{\bs\alpha}'$ for all $\bs\alpha\in\Bbb Z^{\oplus \Bbb N}$ and 
$\dim_F W_{\bs\alpha}/W_{\bs\alpha}'=\dim_F Y_{\bs\alpha}/Y_{\bs\alpha}'=1$ for all $\bs\alpha\neq\bk
0$.) Thus
$(N,N_{\bk 0},{\Cal B})$ is an admissible triple. Since $Y_{\bk 0}$ is the centre, (5) implies that
$N$ is a null system of null rank $\infty$.

We show that this $N$ is indecomposable. In fact, suppose that $N$ is decomposable. Then
$N=(L_1,H_1,{\Cal B}_1)\oplus (L_2,H_2,{\Cal B}_2)$. If $d_1\in H_1$, then $d_1\in H_2$ since $(d_1,
d_1)=0$, which is impossible. Hence,
$d_1\notin H_1$ and $d_1\notin H_2$. So there exist $h_1:=\sum_i\ (u_id_i+v_i\delta_i)\in H_1$ with
$u_1\neq 0$ and
$h_2:=\sum_i\ (p_id_i+q_i\delta_i)\in H_2$  with $p_1\neq 0$ ($u_i,v_i,p_i,q_i\in F$). Note that
$N_{\delta_1}=Ft_1d+Ft_1c$. For $0\neq at_1d+bt_1c\in N_{\delta_1}$ ($a,b\in F$), we have
$[h_1,at_1d+bt_1c]=u_1(at_1d+bt_1c)\neq 0$ since $\delta_i$ are central. Hence,
$at_1d+bt_1c\notin L_2$. Similarly,
$[h_2,at_1d+bt_1c]=p_1(at_1d+bt_1c)\neq 0$, and hence $at_1d+bt_1c\notin L_1$.  But we have
$N_{\delta_1}=(N_{\delta_1}\cap L_1)\oplus(N_{\delta_1}\cap L_1)$, which implies
$N_{\delta_1}=0$, a contradiction.
\endexample

\example{Remark 7.4} In the notation from \S 5, let $M=\frak g\otimes_FF[t_i^{\pm 1}]_{i\in\Bbb N}$.
Then
$M\oplus N$ becomes a LEALA with suitable multiplication. ($H=\frak h\oplus N_\bk 0$  is
infinite-dimensional.) Also, the subalgebra $M\oplus Y'$ is a central extension of $M$. More precisely,
$Y'$ is the centre of $M\oplus Y'$, and the multiplication is given by
$$ [x\otimes t_{\bs\alpha},y\otimes t_{\bs\beta}] =[x,y]\otimes t_{\bs\alpha+\bs\beta}+
\kappa(x,y)t_{\bs\alpha+\bs\beta} c_{\bs\beta}
$$ for $x\otimes t_{\bs\alpha},y\otimes t_{\bs\beta}\in M$. The bracket of $W'$ and $Y'$ is the same as
in $N$, and the bracket of $W'$ and $M$  is just the action of $W'$ on
$F[t_i^{\pm 1}]_{i\in\Bbb N}$ as derivations, i.e.,
$$ [t_{\bs\alpha}d, x\otimes t_{\bs\beta}] =x\otimes d(\bs\beta)t_{\bs\alpha+\bs\beta}
$$ for 
$t_{\bs\alpha}d\in W'$ and
$x\otimes t_{\bs\beta}\in M$. 

Note that this $M\oplus N$ corresponds to the trivial cocycle from $W'\times W'\longrightarrow Y'$ or
in other words the trivial abelian extension of $W'$ by $Y'$. There exists a nontrivial extension
$N_\tau=W'\oplus Y'$ determined by the Moody-Rao cocycle $\tau$ [EM] so that
$N_\tau$ is a null system and
$M\oplus N_\tau$ is a LEALA. More precisely,
$\tau:W'\times W'\longrightarrow Y'$ is defined by
$$ \tau(t_{\bs\alpha}d, t_{\bs\beta}d') 
=d'(\bs\alpha)d(\bs\beta)t_{\bs\alpha+\bs\beta}c_{\bs\alpha-\bs\beta},
$$ and the new bracket
$[\cdot,\cdot]_\tau$ of $N_\tau$ is defined as
$$ [t_{\bs\alpha}d, t_{\bs\beta}d']_\tau =  [t_{\bs\alpha}d, t_{\bs\beta}d'] +\tau(t_{\bs\alpha}d,
t_{\bs\beta}d'),
\tag$1_\tau$ $$ and the rest of brackets remain same as in $N$. As stated in [BGK, Rem.3.76], the
classification of such cocycles seems to be interesting (see also [BeB]).
\endexample

\head
\S 8 Classification for LEALAs of nullity 0
\endhead

Using the results in [NS] and [Stu], we can classify LEALAs of nullity 0. We define the core in
general, and then classify the cores of nullity 0  as the first step.

\definition{Definition 8.1} Let $L$ be a LEALA. The {\it core} of $L$, denoted by $L_c$, is the
subalgebra of $L$ generated by the root spaces $L_\alpha$ for all $\alpha\in R^\times$.
\enddefinition

\example{Remark 8.2}  One can show that the core $L_c$ of a LEALA $L$ is an ideal of $L$ and the
centralizer $C_L(L_c)$ of $L_c$ in $L$ is orthogonal to $L_c$ relative to ${\Cal B}$ (see [BGK,
Lem.3.6] where they gave a proof for EALAs).
\endexample

As in the original EALAs, the core of a LEALA has a well-behaved grading structure (see [AG]). We use
the terminology,
$V$, $V^0$, $\bar V$, and $\bar R^\times$ defined in \S 3.  Recall that $\bar R^\times$ is a locally
finite irreducible root system.  Note that for each $\bar\alpha\in\bar R^\times$ there exists $v\in
V^0$ such that
$\alpha+v\in R$.  So we pick $\dot\alpha\in\bar\alpha\cap R$ for each $\bar\alpha\in\bar
R^\times_\red$, where
$$\bar R^\times_\red=\{\bar\alpha\in\bar R^\times\mid\frac{1}{2}\bar\alpha\notin\bar
R^\times\}\quad\text{(a reduced locally finite irreducible root system)}$$ and
$\dot\beta\in\bar\beta$ for each $\bar\beta\in\bar R^\times\setminus\bar R^\times_\red$ so that
$\dot\beta=2\dot\alpha$ for some $\dot\alpha\in R$ (which is already chosen form $\bar\alpha\in\bar
R^\times_\red$), and fix them. (Note $\dot\beta=2\dot\alpha\notin R$.)  Let $\dot
R=\{\dot\alpha\mid\bar\alpha\in\bar R^\times\}$ and $\dot R_\red=\{\dot\alpha\mid\bar\alpha\in\bar
R^\times_\red\}$, which are clearly a locally finite irreducible root system and a reduced  locally
finite irreducible root system, respectively.  Then
$$R^\times=\bigcup_{\dot\alpha\in\dot R}\ (\dot\alpha+S_{\dot\alpha}),$$ where $S_{\dot\alpha}=\{v\in
V^0\mid\dot\alpha+v\in R^\times\}$ (cf. Remark 4.3).  In particular, $0\in S_{\dot\alpha}$ for
$\dot\alpha\in\dot R_\red$ and $0\notin S_{\dot\beta}$ for $\dot\beta\in\dot R\setminus\dot R_\red$. 
Also, one can show that
$S_{\dot\alpha}\subset R^0$ for all for $\dot\alpha\in\dot R$ using $sl_2$-theory. (For both cases
$S_{\dot\alpha}$ and $S_{2\dot\alpha}$, take an $sl_2$-triplet of the root $\dot\alpha$, and act it on
$L_{\dot\alpha+v}$ or
$L_{2\dot\alpha+v}$.  In fact one can show several more properties for the sets $S_{\dot\alpha}$'s, for
example,
$S_{2\dot\alpha}\subset S_{\dot\alpha}$.)

Now, let $\Omega=\Omega(\dot R)$ be the additive subgroup of $V$ generated by
$\dot R$ and let
${\Lambda}$ be the additive subgroup of $V$ generated by
$\cup_{\dot\alpha\in\dot R}\ S_{\dot\alpha}$. Let 
$$(L_c)_\omega^\delta:=L_c\cap L_{\omega+\delta}$$ for $\omega\in\Omega$ and $\delta\in {\Lambda}$,
defining
$L_{\omega+\delta}=0$ if $\omega+\delta\notin R$. Then one can show that
$$L_c=\bigoplus_{\omega\in\Omega}\ \bigoplus_{\delta\in {\Lambda}}\ (L_c)_\omega^\delta,$$ which is an
$\Omega\oplus {\Lambda}$-graded Lie algebra satisfying:
\roster
\item 
$\dim_F (L_c)_{\dot\alpha}^\delta\leq 1$ for $\dot\alpha\in\dot R$ with
$\dim_F (L_c)_{\dot\alpha}^0= 1$ if ${\dot\alpha}\in \dot R_\red$.

\item If $\dim_F (L_c)_{\dot\alpha}^\delta= 1$, then there exist
$x\in (L_c)_{\dot\alpha}^\delta$ and
$y\in (L_c)_{-\dot\alpha}^{-\delta}$ such that
$[[x,y],u] = \frac{2(\omega,\dot\alpha)}{(\dot\alpha,\dot\alpha)} u$ for $u\in (L_c)_{\omega}^\lambda$
if
$\omega\in\dot R\cup\{0\}$ and $\lambda\in {\Lambda}$.

\item For $\delta\in {\Lambda}$ we have 
 $(L_c)_0^\delta=\sum_{\dot\alpha\in\dot R, \lambda\in {\Lambda}}\  [(L_c)_{\dot\alpha}^\lambda,
(L_c)_{-\dot\alpha}^{\delta-\lambda}]$.

\endroster ((1) follows from Lemma 2.1 and our definition. For (2), take an $sl_2$-triplet for
$\dot\alpha$, and (3) follows from the definition of the core.) Thus,
$L_c$ is in the class of a natural generalization of Lie ${\Lambda}$-tori defined  by Neher [N2] or
[N3]. (If
${\Lambda}$ is a free abelian group of finite rank and if $\dot R$ is a finite irreducible root system,
then such double-graded Lie algebras are exactly the Lie tori he defined.)   Lie ${\Lambda}$-tori have
also been defined by the second author in [Y2] or [Y3], using the notion of root-graded Lie algebras. 
They are equivalent but the Neher's definition is more convenient for this context.

\proclaim{Proposition 8.3} Let $(L, H, {\Cal B})$ be a LEALA.  Then: 

(1) Any finite subset of the core is contained in a Lie torus.   In particular, the core is a direct
limit of Lie tori.

(2) If $L$ has nullity $0$, then the core $L_c$ is a locally finite split simple Lie algebra, and
$(L_c, H\cap L_c, {\Cal B}\mid_{L_c})$ is also a LEALA (of nullity $0$).


\endproclaim
\demo{Proof} Without loss of generality, we may assume that the given subset $S$ consists of homogenous
elements of degree $(\dot\alpha,\delta)$ for $\dot\alpha\in\dot R$ and $\delta\in {\Lambda}$ in the
$\Omega\oplus {\Lambda}$-graded Lie algebra
$L_c$, say $S=\{x_{\dot\alpha_1}^{\delta_1},\ldots, x_{\dot\alpha_r}^{\delta_r}\}$.  Let $Q$ be a
finite full irreducible subsystem of (the locally finite irreducible root system) $\dot R$ containing
$\{\dot\alpha_1,\ldots,
\dot\alpha_r\}$ (see Lemma 3.9).  Then the subalgebra $M$ of $L_c$ generated by
$L_{\dot\alpha_1}^{\delta_1}$,
$\ldots$, $L_{\dot\alpha_r}^{\delta_r}$ and $L_{\dot\alpha_1}^0$, $\ldots$, $L_{\dot\alpha_r}^0$
($L_{\dot\alpha_i}^0=0$ if $\dot\alpha_i\notin\dot R_\red$) is a Lie ${\Lambda}'$-torus of type $Q$
containing
$S$, where
${\Lambda}'$ is the subgroup of ${\Lambda}$ generated by $\{{\delta_1},\ldots, {\delta_r}\}$, which is
a free abelian group of finite rank.  Hence (1) holds.

If $L$ has nullity 0, then ${\Lambda}={\Lambda}'=0$, i.e., the subalgebra $M$ is a Lie 0-torus, which
is a finite-dimensional split simple Lie algebra (by Serre's Theorem).  Hence $L_c$ is locally finite
and simple.  Also, in this case, $R^\times=\dot R$ and $L_c=H'\oplus(\oplus_{\alpha\in R^\times}\
L_\alpha)$, where $H':=H\cap L_c$.  Hence $L_c$ is split relative to $H'$.  Note that ${\Cal
B}\mid_{L_c}$ is nondegenerate since $L_c$ is simple. Thus the second statement of (2) is clear.
%
\qed
\enddemo

To complete the classification of LEALAs of nullity 0, we need to determine the complement of the core. 
We first prove the following:

\proclaim{Lemma 8.4} Let $(L, H, {\Cal B})$ be a LEALA of nullity 0 and $Z$ the centre of $L$. Then 
$Z$ is split and 
$L=L_c\oplus D\oplus Z$ with
$H=(H\cap L_c)\oplus D\oplus Z$, where 
$D$ is an abelian subalgebra acting on $L_c$. In particular, $L_c\oplus D$ is a centreless Lie algebra.
\endproclaim
\demo{Proof} Note first that the centre $Z$ is always contained in $H$ for any LEALA. Also, we have
$Z\cap L_c=0$ since $L_c$ is simple (by Lemma 8.3). We fix a complement $D$ of $(H\cap L_c)\oplus Z$ in
$H$ as a vector space (which is automatically an abelian subalgebra). Then, 
$L=L_c\oplus D\oplus Z$ since $L$ has nullity 0. Now the rest of statements are clear.
\qed
\enddemo

Thus the problem now is the classification of  centreless Lie algebras
$L_c\oplus D$, where $L_c$ is a locally finite split simple Lie algebra and $D$ is an abelian Lie
algebra consisting of derivations of $L_c$ preserving the root spaces of $L_c$ with $[H\cap L_c,
D]=0$.  If
$L_c$ is finite-dimensional, then any derivation is inner. Hence, $L_c\oplus D$ cannot be centreless
unless $D=0$. Thus $D=0$ in this case.  For the infinite-dimensional case, we prove the following
general lemma.

\proclaim{Lemma 8.5} Let $L=\frak g\oplus D$ be a Lie algebra, where $\frak g=\frak
h\oplus\bigoplus_{\mu\in\Delta}\ \frak g_\mu$ is a locally finite split simple Lie algebra with a split
Cartan subalgebra $\frak h$, and $D$ is an abelian subalgebra acting on $\frak g$ with $[\frak h, D]=0$
and $[d,\frak g_\mu]\subset \frak g_\mu$ for all $d\in D$ and $\mu\in\Delta$.  Let ${\sigma}$ be a
symmetric invariant bilinear form on $L$ so that ${\sigma}\mid_{\frak g\times\frak g}\neq 0$.  Then
$\rad {\sigma}\subset Z$, where $\rad {\sigma}$ is the radical of ${\sigma}$ and $Z$ is the centre of
$L$.  In particular, if $Z=0$, then ${\sigma}$ is nondegenerate.

Also, we have $Z=0$ $\Longleftrightarrow$ $D\cap Z=0$ and $(\ad_\frak g \frak g)\cap(\ad_\frak g D)=0$,
i.e., $D$ consists of outer derivations of $\frak g$.
\endproclaim

\demo{Proof}  Clearly, one can assume that $\rad {\sigma}\subset \frak h\oplus D$.  Let
$h+d\in\rad{\sigma}$ for $h\in\frak h$ and $d\in D$.  Let ${\rho}$ be a nondegenerate symmetric
invariant bilinear form on $\frak g$.  For each $\mu\in\Delta$, 
let  $(x_\mu, t_{\mu}, x_{-\mu})$ be a triplet
so that ${\rho}(x_\mu, x_{-\mu})=1$ (see (0.1)).  Note that a symmetric invariant
bilinear form on $\frak g$ is unique up to scalars (see [NS, Lem.II.11]), and so ${\sigma}\mid_{\frak
g\times\frak g}=u{\rho}$ for some $0\neq u\in F$.  Also, $\ad d(x_\mu)= vx_\mu$ for some $v\in F$. 
Now, we have
$0={\sigma}(h+d,t_\mu)=u{\rho}(h,t_\mu)+{\sigma}(d,[x_\mu,x_{-\mu}])
=u\mu(h)+u{\rho}([d,x_\mu],x_{-\mu})=u\mu(h)+uv{\rho}(x_\mu,x_{-\mu})=u\mu(h)+uv$.  
Hence, $v=-\mu(h)$, and so $\ad d=-\ad h$, that is, $h+d\in Z$.

For the second statement, note that $Z\subset \frak h\oplus D$.  Thus $0\neq h+d\in Z$
$\Leftrightarrow$ $h=0$ and $0\neq d\in Z$ or $h\neq 0$ and $d\neq 0$ with $\ad (h+d)=0$
$\Leftrightarrow$ $D\cap Z\neq 0$ or $(\ad_\frak g \frak g)\cap(\ad_\frak g D)\neq 0$.  \qed
\enddemo

Thus, in our question, $D$ has to be contained in the Lie algebra of outer derivations of $L_c$.

Conversely, for any abelian subalgebra
$D$ of outer derivations of a locally finite split simple Lie algebra $\frak g=\frak
h\oplus\bigoplus_{\mu\in\Delta}\ \frak g_\mu$ preserving the root spaces with
$[\frak h, D]=0$, any nondegenerate symmetric invariant bilinear form $\rho$ on $\frak g$ and any
symmetric bilinear form ${\psi}$ on $D$, one can define a centreless LEALA
$(\frak g\oplus D,\frak h\oplus D, {\Cal B})$, where ${\Cal B}=\rho+{\psi}$ with ${\Cal B}(d,\frak
g_\mu)={\Cal B}(\frak g_\mu,d)=0$  and ${\Cal B}(d,t_\mu)={\Cal B}(t_\mu,d)=v_\mu$ for all $d\in D$ and
$\mu\in\Delta$, and $v_\mu$ is a unique element in $F$ so that $[d,x_\mu]=v_\mu x_\mu$ for all
$x_\mu\in\frak g_\mu$. Indeed, all you need to check is the invariance of ${\Cal B}$,  which is easy,
e.g., ${\Cal B}([d,x_\mu],  x_{-\mu})={\Cal B}(v_\mu x_\mu,x_{-\mu})=v_\mu\rho(x_\mu,x_{-\mu})=v_\mu$,
on the other hand, ${\Cal B}(d,[x_\mu,x_{-\mu}])={\Cal B}(d,t_{\mu})=v_\mu$, 
choosing a triplet $(x_\mu, t_{\mu}, x_{-\mu})$ so that $\rho(x_\mu,x_{-\mu})=1$
(see (0.1)).  
Or ${\Cal B}([x_\mu,d], x_{-\mu})={\Cal B}(x_\mu,[d,  x_{-\mu}])$ 
since $[d,x_{-\mu}]=-v_\mu x_\mu$.

Now, assuming that the core $L_c$ is infinite-dimensional, by [NS, Thm IV.6 and Cor.IV.5], 
$(L_c,L_c\cap H)$ is isomorphic to $(\frak g_{\frak I},\frak h)$, where $\frak g_{\frak I}$ is one of 
$sl_{\frak I}(F)$, $o_{2{\frak I}+1}(F)$, $sp_{2{\frak I}}(F)$ or $o_{2{\frak I}}(F)$ for an infinite
index set ${\frak I}$ as defined in Definition 6.1, and $\frak h$ is its standard Cartan subalgebra,
that is,  the subalgebra consisting of diagonal matrices. (They showed that $o_{2{\frak I}+1}(F)\cong
o_{2{\frak I}}(F)$ but $(o_{2{\frak I}+1}(F), \frak h)\not\cong (o_{2{\frak I}}(F), \frak h)$.)  We
identify them, and so $D$ is an abelian subalgebra consisting of outer derivations of $\frak g_{\frak
I}$ preserving the root spaces of $\frak g_{\frak I}$ with $[\frak h, D]=0$.  

By the communication with Karl-Hermann Neeb, we can describe our $D$ more concretely.  For each $d\in
D$, consider
$\tau(x,y):={\Cal B}([d,x],  y)=\rho([d,x],  y)$ for $x,y\in\frak g_{\frak I}$, which is a 2-cocycle of
$\frak g_{\frak I}$ into the trivial module $F$.   Note that the second cohomology
$H^2(\frak g_{\frak I}, F)$ has to be 0 since $\frak g_{\frak I}$ is a direct limit of
finite-dimensional split simple Lie algebras.  Thus there exists a linear form $f$ of $\frak g_{\frak
I}$ such that $\tau(x,y)=f([x,y])$.  Let $A=(a_{ij})$ be a unique matrix 
(of size $\frak I$, $2\frak I+1$ or $2\frak I$) so that $f(x)=\tr(Ax)$ for all
$x\in\frak g_{\frak I}$, where $\tr$ is the trace of the matrix $Ax$.   (Note that only finitely many
diagonal entries of $Ax$ are 0, and so the trace is defined.)  As mentioned in \S 6,  the trace form
$\tr(x,y):=\tr(xy)$ is a symmetric nondegenerate invariant form on $\frak g_{\frak I}$.  Hence, by [NS,
Lem.II.11], $\rho=u\tr$ for some $0\neq u\in F$.  Then we have $$u\tr([d,x], y)={\Cal B}([d,x], 
y)=\tau(x,y)=f([x,y])=\tr(A[x,y])=\tr([A, x]y)$$  for all $x,y\in\frak g_{\frak I}$.   Hence $\ad
d=u\ad A$ if $[A, x]\in\frak g_{\frak I}$. First, we note that $$\text{ if $\frak g_{\frak I}\neq
sl_{\frak I}(F)$, then $a_{ii}=-a_{{\frak I}+i,{\frak I}+i}$, and $a_{2{\frak I}+1,2{\frak I}+1}=0$ for
$o_{2{\frak I}+1}(F)$. }\tag8.6
$$ Thus, for any $\frak g_{\frak I}$ and distinct $i\neq j$,  there exists $k$ distinct from $i,j$ such
that $e_{jj}-e_{kk}\in\frak h$. ($e_{ij}$'s are matrix units.) So we have
$$a_{ij}=\tr(Ae_{ji})=\tr(A[e_{jj}-e_{kk}, e_{ii}]) =u\tr([d,e_{jj}-e_{kk}]e_{ii})=0$$  since $[D,
\frak h]=0$. Therefore, $A$ is diagonal, and in particular, $[A, x]\in\frak g_{\frak I}$.

Thus we have shown that
$D$ consists of diagonal matrices satisfying (8.6) since our Lie algebra is centreless.  Also,
infinitely many
$a_{ii}$ are nonzero since $d$ is outer.  For the case $\frak g_{\frak I}= sl_{\frak I}(F)$, if the
diagonal matrix $A$ is almost scalar, i.e., $a_{ii}=a$ for some $a\in F$ except for finitely many $i$, 
then $\ad A=\ad B$ for $B:=A-\sum_{i\in\frak I}\ ae_{ii}$, and $\ad B=\ad B'$ for some matrix
$B'$ with $\tr(B')=0$. So $\ad A$ is inner on $\frak g_{\frak I}$.  Thus $A$ cannot be an almost scalar
matrix.

Conversely, if $\frak g_{\frak I}\neq sl_{\frak I}(F)$, then the adjoint of any diagonal matrix
satisfying (8.6) with infinitely many nonzero entries is outer, and if $\frak g_{\frak I}= sl_{\frak
I}(F)$, then  the adjoint of any diagonal matrix which is not almost scalar is outer.

Thus we obtain a complete classification of LEALAs of nullity 0.

\proclaim{Theorem 8.7} Let $(L, H, {\Cal B})$ be a LEALA of nullity $0$ and $Z$ the centre of $L$.

(1) If $L_c$ is finite-dimensional, then 
$(L, H, {\Cal B})\cong (\frak g\oplus Z,\frak h\oplus Z,u\kappa\perp {\psi})$ for some nonzero
$u\in F$,  where
$\frak g$ is a finite-dimensional split simple Lie algebra,
$\frak h$ is a split Cartan subalgebra,
$\kappa$ is the Killing form of $\frak g$, and ${\psi}$ is a nondegenerate symmetric bilinear form on
$Z$. In particular, $L$ is a split central extension of $\frak g$.

Conversely, any such $Z$, $u$ and ${\psi}$ give a LEALA 
$(\frak g\oplus Z,\frak h\oplus Z,u\kappa\perp {\psi})$ (an EALA if $\dim_F Z<\infty$).

\enskip

(2) If $L_c$ is infinite-dimensional, then 
$(L, H, {\Cal B})\cong (\frak g_{\frak I}\oplus D\oplus Z,\frak h\oplus D\oplus Z,u\tr+{\psi})$ for
some nonzero $u\in F$,  where
$\frak g_{\frak I}=sl_{\frak I}(F)$, $o_{2{\frak I}+1}(F)$, $sp_{2{\frak I}}(F)$ or $o_{2{\frak I}}(F)$
with its standard Cartan subalgebra $\frak h$ for an infinite index set ${\frak I}$, and
$D$ is:
\roster
\item"(i)" a subspace consisting of  diagonal matrices  satisfying (8.6) with infinitely many nonzero
entries if 
$\frak g_{\frak I}\neq sl_{\frak I}(F)$; 
\item"(ii)"  a subspace consisting of  diagonal matrices which are not almost scalar if 
 $\frak g_{\frak I}=sl_{\frak I}(F)$,
\endroster and ${\psi}$ is a symmetric bilinear form on
$D\oplus Z$ satisfying $Z\cap \rad {\psi}=0$. In particular,
$L$ is a split central extension of the centreless LEALA
$(\frak g_{\frak I}\oplus D, \frak h\oplus D, u\tr+{\psi}\mid_{ D\times
 D})$.

Conversely, any such $Z$, $ D$, $u$ and ${\psi}$ give a LEALA 
$(\frak g_{\frak I}\oplus  D\oplus Z,\frak h\oplus  D\oplus Z,u\tr+{\psi})$.

\endproclaim
\demo{Proof}  We only need to prove that $u\tr+{\frak c}$ is nondegenerate if and only if $Z\cap \rad
{\psi}=0$. Clearly, if $Z\cap \rad {\psi}\neq 0$, then $u\tr+{\psi}$ is degenerate. Thus we need to
prove that $Z\cap \rad {\psi}=0$ implies the nondegeneracy of $u\tr+{\psi}$. Let
${R}$ be the radical of $u\tr+{\psi}$,  which is clearly contained in $ D\oplus Z$. Let $A+z\in {R}$
for $A\in D$ and $z\in Z$. If $A\neq 0$, then there exists $0\neq x\in\frak g_{\frak I}$ such that
$[x,A+z]=vx\in{R}$  for some $0\neq v\in F$ since the radical is an ideal. This implies that $\frak
g_{\frak I}\subset {R}$, which is a contradiction.  Hence, $A=0$, and so ${R}\subset Z$, which implies
that ${R}\subset Z\cap \rad {\psi}$. Thus we get ${R}=0$.
\qed
\enddemo

Note that we gave some of the smallest nontrivial examples in Example 6.2.

\head
\S 9 A note for tame LEALAs
\endhead

We can define the tameness for a LEALA as in [A-P] (or [N2]).
\definition{Definition 9.1} Let $L$ be a LEALA. Let  $C=C_L(L_c)$ be  the centralizer. Then
$L$ is called  {\it tame} if $C\subset L_c$.
\enddefinition 

By Theorem 8.7, if the nullity is 0, then the centre has to be 0. So:
\proclaim{Corollary 9.2} Let $L$ be a tame LEALA of nullity $0$. Then $L\cong\frak g$ or 
$\frak g_{\frak I}\oplus D$ in the description of Theorem 8.7.  In particular,  tame
$\Longleftrightarrow$ centreless
$\Longleftrightarrow$ indecomposable, in the case of nullity $0$. 
\endproclaim However, this is not the case if the nullity is bigger than 0. There are examples of
LEALAs which are indecomposable but not tame in [BGK]. For the convenience of the reader, we give one
such example.

\example{Example 9.3} Let $M=\frak g\otimes F[t^{\pm 1}]$ be a loop algebra with form
$\kappa\otimes\varepsilon$ (see \S 5) and
$A=M\oplus Fc\oplus Fd$ its (untwisted) affine Kac-Moody Lie algebra (cf. Remark 7.2). Let
$(N,N_0,{\Cal B})=\oplus_{\delta\in S}\ N_\delta$ be a null system constructed from 
${\Lambda}=\Bbb Z$, $S$ and $\vvp=\id$ in Example 7.1. We identify $N_0$  with the subalgebra $Fc\oplus
Fd$ of $A$. Let $L$ be the Lie algebra containing $A$ and $N$ as subalgebras so that $L=A\oplus
(\oplus_{\delta\neq 0}\ N_\delta)$ declaring
$[M, \oplus_{\delta\neq 0}\ N_\delta]=0$. Let $H:=(\frak h\otimes 1)\oplus N_0$. Then,
$(L,H,\kappa\otimes\varepsilon+{\Cal B})$ is an EALA of null rank 1. Since an affine Kac-Moody Lie
algebra is indecomposable,
$H$ is never decomposed into a direct sum of orthogonal subspaces. Hence $L$ is an indecomposable EALA.
However, the subspace $\oplus_{\delta\neq 0}\ N_\delta$ centralizes the core of $L$, i.e., $L_c=M\oplus
Fc$. So $L$ is not tame unless $S=0$.
\endexample

Thus, to classify even EALAs of null rank 1, it seems that we need at least to classify null systems of
null rank 1. So it may be natural to assume the tameness. Also, by the following lemma, the notion of
indecomposability disappears if the tameness is assumed.

\proclaim{Lemma 9.4} A tame LEALA is indecomposable. Also, an EALA is completely decomposable with the
factors of one indecomposable EALA  and some indecomposable null systems.
\endproclaim
\demo{Proof} Let  $L=(L,H,{\Cal B})$ be a LEALA. Suppose that 
$L=(L_1,H_1,{\Cal B}_1)\oplus (L_2,H_2,{\Cal B}_2)$. Then 
$R(H)=R(H_1)\cup R(H_2)$ and $(R(H_1), R(H_2))=0$. So by (A5),
$R^\times\subset R(H_1)$ or $R(H_2)$. Thus the core $L_c$ has to sit in one of the factors, say
$L_c\subset L_1$. Then $L_2$ is a null system. So if $L$ is tame, then $L_2=0$ since $L_2$ centralizes
$L_c$.

For the second statement,  if $L$ is an EALA, then $L$ is completely decomposable since $\dim
H<\infty$. Thus, by the same reason as above, the core $L_c$ has to be in one of the factors, and the
rest are null systems.
\qed
\enddemo


\Refs \widestnumber\key{ABGP}

\ref
\key A-P
\by B.N.  Allison, S.  Azam, S. Berman, Y.  Gao, A. Pianzola
\book Extended Affine Lie Algebras and Their Root Systems
\bookinfo Memoirs Amer.  Math.  Soc.
\text{\bf 126}
\vol 603
\yr 1997
\endref



\ref 
\key AG
\by B.N. Allison, Y. Gao 
\paper The root system and the core of an extended affine Lie algebra
\jour Selecta Mathematica, New Series 
\vol 7 
\yr 2001 
\pages 149--212
\endref



\ref
\key AKY
\by S. Azam, V. Khalili, M. Yousofzadeh
\paper Extended affine root systems of type BC
\jour preprint
\vol 
\yr 
\pages 
\endref


\ref
\key AY
\by B.N. Allison, Y. Yoshii
\paper Structurable tori and extended affine Lie algebras of type 
$\text{BC}_1$
\jour  J. Pure Appl. Algebra
\vol 184(2-3)
\yr 2003
\pages 105--138
\endref




\ref\key B
\by R.E. Block
\paper  The irreducible representations  of the Lie algebra $\frak s\frak l(2)$ and of the Weyl algebra
\jour Adv. in Math. 
\vol 39(1)
\yr 1981
\pages 69--110
\endref


\ref
\key Bo
\by N. Bourbaki
\book Groupes et alg\` ebres de Lie, Chap. IV, V, VI
\bookinfo 
\publ Hermann, Paris
\yr 1968
\endref




\ref\key BGK
\by S. Berman, Y. Gao, Y. Krylyuk 
\paper  Quantum tori and the structure of elliptic quasi-simple Lie algebras
\jour J. Funct. Anal. 
\vol 135
\yr 1996
\pages 339--389
\endref

\ref
\key BB
\by Y. Bahturin, G. Benkart
\paper Some constructions in the theory of locally finite simple Lie algebras
\jour  J. Lie Theory
\vol 14
\yr 2004
\pages 243--270
\endref

\ref
\key BeB
\by S. Berman, Y. Billig
\paper Irreducible representations for toroidal Lie algebras
\jour  J. Algebra
\vol 221(1)
\yr 1999
\pages 188--231
\endref


\ref\key HT
\by R. H\o egh-Krohn, B. Torresani
\paper  Classification and construction of quasi-simple Lie algebras
\jour J. Funct. Anal.
\vol 89
\yr 1990
\pages 106--136
\endref

\ref
\key EM
\by R. Eswara, R. Moody
\paper Vertex representations for $n$-toroidal Lie algebras and a generalization of the Virasoro algebra
\jour Comm. Math. Phys.
\vol 159(2)
\yr 1994
\pages 239--264
\endref 



\ref
\key G
\by Y.  Gao
\paper The degeneracy of extended affine Lie algebras
\jour
 Manuscripta Math.
\vol 97(2)
\yr 1998
\pages 233--249
\endref 


\ref
\key GN
\by E. Garcia, E. Neher
\paper Gelfand-Kirillov  dimension and local finiteness of Jordan superpairs  covered by grids and
their associated Lie superalgebas
\jour Comm. Alg.
\vol 32
\yr 2004
\pages 2149--2176
\endref 




\ref
\key K
\by V. Kac
\book Infinite dimensional Lie algebras
\bookinfo third edition
\publ Cambridge University Press
\yr 1990
\endref


\ref
\key LN
\by O. Loos, E. Neher
\book Locally finite root systems
\bookinfo Memoirs Amer.  Math.  Soc.
\text{\bf 811}
\vol 171
\yr 2004
\endref 


\ref
\key MP
\by R.V. Moody, A. Pianzola
\book Lie algebras with triangular decompositions
\bookinfo Can. Math. Soc. series of monographs and advanced texts
\publ John Wiley
\yr 1995
\endref


\ref 
\key N1
\by E. Neher
\paper  Lie algebras graded by 3-graded root systems
\jour Amer. J. Math.
\vol 118
\yr 1996
\pages 439--491
\endref

\ref
\key N2
\by E.  Neher
\paper Extended affine Lie algebras
\jour C. R. Math. Rep. Acad. Sci. Canada
\vol 26(3)
\yr 2004
\pages 90--96
\endref 

\ref
\key N3
\by E.  Neher
\paper Lie tori
\jour C. R. Math. Rep. Acad. Sci. Canada
\vol 26(3)
\yr 2004
\pages 84--89
\endref 

\ref
\key NP
\by K.-H. Neeb, I. Penkov
\paper Cartan subalgebras of $gl_\infty$
\jour
 Canad. Math. Bull.
\vol 46(4)
\yr 2003
\pages 597--616
\endref 



\ref
\key NS
\by K.-H. Neeb, N. Stumme
\paper The classification of locally finite split simple Lie algebras
\jour
 J. reine angew. Math.
\vol 533
\yr 2001
\pages 25--53
\endref 


\ref
\key S
\by K. Saito
\paper Extended affine root systems 1 (Coxeter transformations)
\jour RIMS., Kyoto Univ.
\vol 21
\yr 1985
\pages 75--179
\endref


\ref
\key St
\by H. Strade
\paper Locally finite-dimensional Lie algebras and their derivation algebras
\jour Abh. Math. Sem. Univ. Hamburg
\vol 69
\yr 1999
\pages 373--391
\endref 



\ref
\key Stu
\by N. Stumme
\paper On the structure of locally finite split Lie algebras
\jour
 J. Algebra
\vol 220
\yr 1999
\pages 664--693
\endref 


\ref
\key Y1
\by Y. Yoshii
\paper Coordinate algebras of extended affine Lie algebras of type $\text A_1$
\jour J. Algebra
\vol 234
\yr 2000
\pages 128--168
\endref
  
\ref
\key Y2
\by Y. Yoshii
\paper Root systems extended by an abelian group and their Lie algebras
\jour J. Lie Theory
\vol 14(2)
\yr 2004
\pages 371--394
\endref

\ref
\key Y3
\by Y. Yoshii
\paper Lie tori -- A simple characterization of  extended affine Lie algebras
\jour submitted
\endref

\endRefs
\enddocument